\documentclass{icmart}

\contact[entov@math.technion.ac.il]{Department of Mathematics\\
Technion -- Israel Institute of Technology\\
Haifa 32000, Israel}

\newtheorem{thm}{Theorem}[section]
\newtheorem{cor}[thm]{Corollary}

\theoremstyle{definition}

\newtheorem{rem}[thm]{Remark}
\newtheorem{question}[thm]{Question}
\newtheorem{exam}[thm]{Example}

\newcommand{\Ker}{\operatorname{Ker}}

\newcommand{\R}{{\mathbb{R}}}

\newcommand{\Q}{{\mathbb{Q}}}

\newcommand{\Z}{{\mathbb{Z}}}
\newcommand{\N}{{\mathbb{N}}}
\newcommand{\C}{{\mathbb{C}}}

\newcommand{\cA}{{\mathcal{A}}}

\newcommand{\cE}{{\mathcal{E}}}
\newcommand{\cF}{{\mathcal{F}}}

\newcommand{\cK}{{\mathcal{K}}}

\newcommand{\cU}{{\mathcal{U}}}

\newcommand{\gog}{{\mathfrak{g}}}
\newcommand{\gou}{{\mathfrak{u}}}
\newcommand{\spn}{{\mathfrak{sp}\, (2n,\R)}}

\newcommand{\tG}{{\widetilde{G}}}

\newcommand{\tphi}{{\widetilde{\phi}}}
\newcommand{\tpsi}{{\widetilde{\psi}}}

\newcommand{\Cal}{{\hbox{\it Cal}}}
\newcommand{\tCal}{{\widetilde{\hbox{\it Cal}}}}

\newcommand{\Ham}{{\hbox{\it Ham\,}}}

\newcommand{\tHam}{\widetilde{\hbox{\it Ham}\, }}

\newcommand{\Symp}{{\hbox{\it Symp}\, }}
\newcommand{\supp}{{\rm supp\,}}

\title[Quasi-morphisms and quasi-states in symplectic topology]{Quasi-morphisms and quasi-states in symplectic topology}

\author[Michael Entov]
{Michael Entov\thanks{Partially supported by the Israel Science
Foundation grant $\#$ 723/10.}}

\begin{document}

\begin{abstract}
We discuss certain ``almost homomorphisms" and ``almost linear" functionals
that have appeared in symplectic topology
and their applications concerning Hamiltonian dynamics,
functional-theoretic properties of Poisson brackets and algebraic
and metric properties of symplectomorphism groups.

\end{abstract}



\maketitle

\section{Introduction}

Symplectic manifolds carry several interesting mathematical
structures of different flavors, coming from algebra, geometry,
topology, dynamics and analysis. In this survey we discuss certain
``almost homomorphisms" (called {\it Calabi quasi-morphisms}) on
groups of symplectomorphisms of symplectic manifolds and certain
``almost linear" functionals (called {\it symplectic quasi-states})
on the spaces of smooth functions on symplectic manifolds that have
been useful in finding new relations between these structures. In
particular, we describe applications of these new tools to
Hamiltonian dynamics, functional-theoretic properties of Poisson
brackets as well as algebraic and metric properties of the groups of
symplectomorphisms. We also briefly discuss a relation between
the symplectic quasi-states and von Neumann's mathematical
foundations of quantum mechanics. We end the survey with a discussion
on the function theory approach to symplectic topology.

A detailed introduction to the subject can be found in the
forthcoming book \cite{PR-book} by L.Polterovich and D.Rosen.

\bigskip
\noindent {\bf Acknowledgements.} Most of the material presented
below is based on our joint papers with Leonid Polterovich -- I
express my deep gratitude to Leonid for the long and enjoyable
collaboration. Some of the results were obtained jointly with Paul
Biran, Lev Buhovsky, Frol Zapolsky, Pierre Py and Daniel Rosen - I
thank them all.

\section{Quasi-morphisms and quasi-states - generalities}

\subsection{Quasi-morphisms}
 \label{subsec-quasi-morphisms}
Let $G$ be a group. A function $\mu: G\to\R$ is called a {\it
quasi-morphism}, if there exists a constant $C>0$ so that $|\mu
(xy)-\mu(x)-\mu(y)|\leq C$ for any $x,y\in G$. We say that a
quasi-morphism $\mu: G\to\R$ is {\it homogeneous}\footnote{Sometimes
homogeneous quasi-morphisms are also called {\it
pseudo-characters} -- see e.g. \cite{Shtern}.}, if $\mu(x^k)=k\mu(x)$ for any $x\in G$, $k\in
\Z$.

Clearly, any $\R$-valued homomorphism on $G$ is a homogeneous
quasi-morphism but finding homogeneous quasi-morphisms that are not
homomorphisms is usually a non-trivial task.
Let us also note that any
homogeneous quasi-morphism $\mu$ is conjugacy-invariant and
satisfies $\mu (xy)= \mu(x) + \mu(y)$ for any commuting elements $x,
y$ (in particular, any homogeneous quasi-morphism on an abelian
group is a homomorphism). For more on quasi-morphisms see e.g. \cite{Calegary-scl}.

\subsection{Quasi-states and quantum mechanics}
\label{subsec-quantum-mechanics}

Roughly speaking, qua\-si-states are ``almost linear" functionals on
algebras of a certain kind. The term ``quasi-state" comes from the work of
Aarnes (see \cite{Aar91} and the references to the Aarnes' previous work therein)
but its history goes back to the mathematical model of quantum
mechanics suggested by von Neumann \cite{von Neumann}. A basic
object of this model is a real Lie algebra of observables that will be denoted by  $\cA_q$
($q$ for quantum): its elements (in the simplest version of the
theory) are Hermitian operators on a finite-dimensional complex
Hilbert space $H$ and the Lie bracket is given by $[A,B]_{\hbar} =
\frac{i}{\hbar}(AB-BA)$, where $\hbar$ is the Planck constant.
Observables represent physical quantities such as energy, position,
momentum etc. In von Neumann's model a state of a quantum system is
given by a functional $\zeta: \cA_q \to \R$ which satisfies the
following axioms:

\medskip
\noindent {\bf Additivity:} $\zeta(A+B) = \zeta(A)+\zeta(B)$ for
all $A,B \in \cA_q$.

\medskip
\noindent {\bf Homogeneity:} $\zeta(cA) = c\zeta(A)$ for all $c \in
\R$ and $A \in \cA_q$.

\medskip
\noindent{\bf Positivity:} $\zeta(A)\geq 0$ provided $A \geq 0$.

\medskip
\noindent {\bf Normalization:} $\zeta(Id) = 1$.

\medskip

As a consequence of these axioms von Neumann proved that for every
state $\zeta$ there exists a non-negative Hermitian operator
$U_\zeta$ with trace $1$ so that $\zeta (A) = \text{tr}(U_{\zeta}A)$
for all $A \in \cA_q$. An easy consequence of this formula is that
for every state $\zeta$ there exists an observable $A$ such that
\begin{equation} \label{eq-disp}  \zeta(A^2) -(\zeta(A))^2 > 0\;.
\end{equation}
In his book \cite{von Neumann} von Neumann adopted a statistical
interpretation of quantum mechanics according to which the value
$\zeta(A)$ is considered as the expectation of a physical quantity
represented by $A$ in the state $\zeta$. In this interpretation the
equation \eqref{eq-disp} says that there are no dispersion-free
states. This result led von Neumann to a conclusion which can be
roughly described as the impossibility to present random
quantum-mechanical phenomena as an observable part of some ``hidden"
underlying deterministic mechanism.
This conclusion
caused a major discussion among physicists (see e.g. \cite{Bell})
some of whom disagreed with the additivity axiom of a quantum state.
Their reasoning was that the formula $\zeta(A+B) =
\zeta(A)+\zeta(B)$ makes sense {\it a priori} only if observables
$A$ and $B$ are simultaneously measurable, that is, commute:
$[A,B]_{\hbar} = 0$.

In 1957 Gleason \cite{Gleason} proved his famous theorem which can
be viewed as an additional argument in favor of von Neumann's
additivity axiom. Recall that two Hermitian operators on a
finite-dimensional Hilbert space commute if and only if they can be
written as polynomials of the same Hermitian operator. Let us define
a {\it quasi-state} on $\cA_q$ as an $\R$-valued functional which
satisfies the homogeneity, positivity and normalization axioms
above, while the additivity axiom is replaced by one of the two {\it
equivalent} axioms:

\noindent {\bf Quasi-additivity-I:} $\zeta (A+B) =
\zeta(A)+\zeta(B)$, provided $A$ and $B$ commute: $[A,B]_{\hbar} =
0$;

\noindent {\bf Quasi-additivity-II:} $\zeta (A+B) =
\zeta(A)+\zeta(B)$, provided $A$ and $B$ belong to a singly generated
subalgebra of $\cA_q$.

According to the Gleason theorem, every quasi-state $\zeta$ on $\cA_q$ is linear, that is, a state, provided the complex dimension of the Hilbert space $H$ is at least 3 (it is an easy exercise to show that in the two-dimensional case there are plenty of non-linear quasi-states).

Let us turn now to the mathematical model of classical mechanics.
Here the algebra $\cA_c$ of observables ($c$ for classical) is the
space $C^\infty (M)$ of smooth functions on a symplectic manifold
$M$. The space $C^\infty (M)$ carries two structures. On one hand,
it is a Lie algebra with respect to the Poisson bracket (see
Section~\ref{subsec-sympl-prelim}). On the other hand, it is a dense
subset (in the uniform norm) of the commutative Banach algebra
$C(M)$ of continuous functions on $M$. For both frameworks one can
define its own version of the notion of a quasi-state adapting,
respectively, the first or the second definition of quasi-additivity
-- as a result one gets the so-called {\it Lie quasi-states} and
{\it topological quasi-states} (see Section~\ref{subsec-Lie-qs}).

{\it Symplectic quasi-states} that appear in symplectic topology and
will be discussed below in Section~\ref{sec-Calabi-qmms-sympl-qst}
belong to both of these worlds -- they are simultaneously Lie and
topological quasi-states\footnote{Interestingly, symplectic
quasi-states had appeared in an infinite-dimensional setting in
symplectic topology before Lie quasi-states were properly studied in
the finite-dimensional setting.}. Note that for the Lie algebra
$C^\infty (M)$ the first definition of quasi-additivity fits in with
the physical Correspondence Principle according to which the bracket
$[\;,\;]_{\hbar}$ corresponds to the Poisson bracket $\{\;,\,\}$ in
the classical limit $\hbar \to 0$. The existence of non-linear
symplectic quasi-states on certain symplectic manifolds (see
Section~\ref{sec-Calabi-qmms-sympl-qst}) can be viewed as an
``anti-Gleason phenomenon" in classical mechanics. Interestingly, at
least for $M=S^2$, the symplectic quasi-state that we construct is
dispersion-free (see Example~\ref{exam-median-qstate}), unlike
states in von Neumann's model of quantum mechanics. For more
information on the connection of symplectic quasi-states to physics
see \cite{EPZ-physics} and Remark~\ref{rem-from-symplectic-to-quantum-mechanics}
below.


\subsection{Lie and topological quasi-states} \label{subsec-Lie-qs} Here
is the precise definition of a Lie quasi-state. Let $\gog$ be a
(possibly infinite-dimensional) Lie algebra over $\R$ and let
$W\subset \gog$ be a vector subspace. A function $\zeta: W \to \R$
is called a {\it Lie quasi-state}, if it is linear on every abelian
subalgebra of $\gog$ contained in $W$.

Finding non-linear Lie quasi-states is, in general, a non-trivial
task: for instance, the difficult Gleason theorem mentioned above is
essentially equivalent -- in the finite-dimensional setting -- to
the claim that any Lie quasi-state on the unitary Lie algebra
$\gou(n)$, $n\geq 3$, which is bounded on a neighborhood of zero,
has to be linear \cite{EP-Lie-qs}. Choosing an appropriate
regularity class of Lie quasi-states is essential for this kind of
results: if $\gog$ is finite-dimensional, then any Lie quasi-state
on $\gog$ which is differentiable at $0$ is automatically linear
while the space of all, not necessarily continuous, Lie quasi-states
on $\gog$ might be infinite-dimensional \cite{EP-Lie-qs}.

Another source of interest to Lie quasi-states
lies in their connection to quasi-morphisms on Lie groups: if
$\gog$ is the Lie algebra of a Lie group $G$ and $\mu: G\to \R$ is a
homogeneous quasi-morphism continuous on 1-parametric subgroups,
then {\it the derivative of $\mu$}, that is, the composition of $\mu$
with the exponential map, is a Lie quasi-state on $\gog$, invariant
under the adjoint action of $G$ on $\gog$. Symplectic quasi-states
that we will discuss below appear as a particular case of this
construction.

Unfortunately, rather little is known about non-linear (continuous)
Lie quasi-states and their connections to quasi-morphisms in general
-- almost all known facts (in particular, a non-trivial description of
the space of all non-linear continuous Lie quasi-states
on $\spn$, $n\geq 3$) and some
basic open questions on the subject can be found in \cite{EP-Lie-qs}.

Let us now define the notion of a topological quasi-state -- it is
due to Aarnes \cite{Aar91} (who called it just a ``quasi-state").
Let $X$ be a compact Hausdorff topological space and let $C(X)$ be
the space of continuous functions on $X$ equipped with the uniform
norm. For a function $F\in C(X)$ denote by $\cA_F$ the closure in
$C(X)$ of the set of functions of the form $p \circ F$, where $p$ is
a real polynomial. A functional $\zeta:
C(X) \to \R$ is called a {\it topological quasi-state} \cite{Aar91},
if it satisfies the following axioms:

\smallskip \noindent {\bf Quasi-linearity:} $\zeta$ is linear on $
\cA_F$ for every $F \in C(X)$ (in particular, $\zeta$ is
homogeneous).

\smallskip \noindent {\bf Monotonicity:} $\zeta (F) \leq \zeta (G)$
for $F \leq G$.

\smallskip \noindent {\bf Normalization:} $\zeta (1) = 1$.

A linear topological quasi-state is called a {\it state} (similarly to states in von Neumann's model of quantum mechanics -- see Section~\ref{subsec-quantum-mechanics}). The
existence of non-linear topological quasi-states was first proved by
Aarnes \cite{Aar91}.

By the classical Riesz representation theorem, states on $C(X)$ are
in one-to-one correspondence with regular Borel probability measures
on $X$. In \cite{Aar91} Aarnes proved a generalized Riesz
representation theorem that associates to each topological
quasi-state $\zeta$ a {\it quasi-measure}\footnote{Quasi-measures
are sometimes also called {\it topological measures}.} $\tau$
on $X$ which is defined only on sets that are either open or closed
and is finitely additive but not necessarily sub-additive. The
relation between $\zeta$ and $\tau$ extends the relation
between states and measures given by the Riesz representation
theorem. In particular, if $A$ is closed, $\tau (A)$ can be
thought of as the ``value" of $\zeta$ on the (discontinuous)
characteristic function of $A$.

\section{Calabi quasi-morphisms, symplectic qua\-si-states}
\label{sec-Calabi-qmms-sympl-qst}

\subsection{Symplectic preliminaries}
\label{subsec-sympl-prelim} Referring the reader to
\cite{McD-Sal-intro} for the foundations of symplectic geometry we
briefly recall the basic notions needed for the further discussion.

Let $M^{2n}$ be a closed connected manifold equipped with a
symplectic form $\omega$, that is, a closed and non-degenerate
differential 2-form. In terms of classical mechanics, $M$ can be
viewed as the phase space of a mechanical system and smooth
functions on $M$ (possibly depending smoothly on an additional
parameter, viewed as time) are called Hamiltonians. Whenever we
consider a time-dependent Hamiltonian we assume that it is
1-periodic in time, i.e. has the form $F: M\times S^1\to\R$. Set
$F_t := F(\cdot, t)$. The support of $F$ is defined as $\supp F :=
\cup_{t\in S^1} \supp F_t \subset M$. We say that $F$ is {\it
normalized}, if $\int_M F_t\omega^n =0$ for any $t\in S^1$.

We denote by $C(M)$ (respectively, $C^\infty (M)$), the space of
continuous (respectively, smooth) functions on $M$ and by
$\|\cdot\|$ the uniform norm on these spaces: $\| F\|:= \max_M |F|$.

Given a (time-dependent) Hamiltonian $F$, define its {\it
Hamiltonian vector field} $X_{F_t}$ by $\omega (\cdot, X_{F_t}) =
dF_t (\cdot)$. Denote the flow of $X_{F_t}$ by $\phi_F^t$ -- it
preserves $\omega$ and is called the {\it Hamiltonian flow of $F$}.
Symplectomorphisms of $M$ (that is, diffeomorphisms of $M$
preserving $\omega$) that can be included in such a flow are called
{\it Hamiltonian symplectomorphisms} and form a group $\Ham (M)$
which is a subgroup of the identity component $\Symp_0 (M)$ of the
full symplectomorphism group $\Symp (M)$. Its universal cover is
denoted by $\tHam (M)$. We say that $\phi_F:=\phi_F^1$ is {\it the
Hamiltonian symplectomorphism generated by $F$}. The Hamiltonian $F$
also generates an element of $\tHam (M)$ that will be denoted by
$\tphi_F$: it is given by the homotopy class (with the fixed
end-points) of the path $\phi_F^t$, $0\leq t\leq 1$, in $\Ham (M)$.

The space of smooth functions on $M$ will be denoted by $C^\infty
(M)$. Given $F,G\in C^\infty (M)$, define the {\it Poisson bracket}
$\{F,G\}$ by $\{ F,G\} := \omega (X_G,X_F)$. Together with the
Poisson bracket $C^\infty (M)$ becomes a Lie algebra whose center is
$\R$ (the constant functions).

It is instructive to view $\Ham (M)$ and $\tHam (M)$ as
infinite-dimensional Lie groups whose Lie algebra (the algebra of
time-independent Hamiltonian vector fields on $M$) is naturally
isomorphic to $C^\infty (M)/\R$, or to the subalgebra of $C^\infty
(M)$ formed by normalized functions, with the map $F\mapsto \phi_F$
being viewed as the exponential map and the natural action of $\Ham
(M)$ on $C^\infty (M)$ being viewed as the adjoint action.

Similarly to the closed case, for an open symplectic manifold
$(U^{2n},\omega)$ one can define $\Ham (U)$ as the group formed by
Hamiltonian symplectomorphisms generated by (time-dependent)
Hamiltonians supported in $U$ and $\tHam (U)$ as its universal
cover. The group $\tHam (U)$ admits the {\it Calabi homomorphism}
$\tCal_U: \tHam (U)\to\R$ defined by $\tCal_U (\tphi_F):= \int_{S^1}
\int_M F_t\omega^n dt$, where $\supp F \subset U$. If $\omega$ is
exact, $\tCal_U$ descends to a homomorphism $\Cal_U: \Ham
(U)\to\R$. If $U$ is an open subset of $M$, then there are natural
inclusion homomorphisms $\Ham (U)\to \Ham (M)$, $\tHam (U)\to\tHam
(M)$, whose images will be denoted by $G_U$ and
$\tG_U$\footnote{Note that the homomorphism $\tHam
(U)\to\tHam (M)$ does not have to be injective and, accordingly,
$\tG_U$ does not have to be the preimage of $G_U$ under the
universal cover $\tHam (M)\to \Ham (M)$ -- see \cite{Kislev}.}.

Let $U$ be an open subset of $M$. Each $\phi\in \Ham (M)$
(respectively, $\tphi\in\tHam (M)$) can be represented as a product
of elements of the form $\psi\theta\psi^{-1}$ with $\theta$ lying in
$G_U$ (respectively, $\tG_U$). Moreover, assuming that $\tCal_U$
descends to a homomorphism $\Cal_U$ on $\Ham (U)$, one can make sure that each such $\theta$ satisfies $\Cal_U
(\theta)=0$. This follows from Banyaga's fragmentation lemma
\cite{Banyaga}. Denote the minimal number of factors in such a
product by $\|\phi\|_U$ (respectively, $\|\tphi\|_U$), if there is
no condition on $\Cal_U (\theta)$, and $\|\phi\|_{U,0}$, if the
condition $\Cal_U (\theta)=0$ is imposed. All the norms are defined
as $0$ on the identity elements.

Let $T^k$ be a torus. A {\it Hamiltonian $T^k$-action} on $M$ is a
homomorphism $T^k\to\Ham (M)$. (We will always assume that such an action
is effective). In such a case the action of the
$i$-th $S^1$-factor of $T^k=S^1\times\ldots\times S^1$,
$i=1,\ldots,k$, is a Hamiltonian flow generated by a Hamiltonian
$H_i$. The Hamiltonians $H_1,\ldots, H_k$ commute with respect to
the Poisson bracket. The map $\Phi=(H_1,\ldots,H_k): M\to\R^k$ is
called the {\it moment map} of the Hamiltonian $T^k$-action. If all
$H_i$ are normalized, we say that $\Phi$ is the normalized moment
map.

A submanifold $L$ of $(M^{2n},\omega)$ is called {\it Lagrangian}, if
${\rm dim}\, L=n$ and $\omega|_L\equiv 0$.

A (closed) symplectic manifold $(M,\omega)$ admits a preferred class
of almost complex structures compatible in a certain sense with $\omega$. All these almost complex structures
have the same first Chern class $c_1$, called the first Chern class of $M$. A closed symplectic manifold
$(M,\omega)$ is called {\it
monotone}, if $[\omega]$ and $c_1$ are positively proportional on
spherical homology classes and {\it symplectically aspherical}, if
$[\omega]$ vanishes on such classes.

Finally, we say that a subset $X\subset M$ is {\it displaceable from
$Y\subset M$ by a group $G$} (where $G$ is either $\Ham (M)$,
or $\Symp_0 (M)$, or $\Symp (M)$), if there exists $\phi\in G$ such
that $\phi (X)\cap \overline{Y} = \emptyset$. If $X$ can be
displaced from itself by $G$, we say that it is {\it
displaceable by $G$} or just {\it displaceable}, if $G=\Ham (M)$.

Consider $T^* S^1 =
\R\times S^1$ with the coordinates $(r, \theta)$ and the
symplectic form $dr \wedge d\theta$. We say that $X\subset M$ is
{\it stably displaceable}, if $X \times \{ r = 0\}$ is
displaceable in $M \times T^* S^1$ equipped with the split
symplectic form $\omega \oplus (dr \wedge d\theta)$. Any displaceable
set is stably displaceable (but not necessarily vice versa).

\subsection{Quantum homology and spectral numbers}
\label{subsec-quantum-hom-spec-numbers} A closed symplectic manifold
$M$ carries a rich algebraic structure called the {\it quantum
homology} of $M$: additively it is just the singular homology of $M$
with coefficients in a certain ring, while multiplicatively the quantum product is a
deformation of the classical intersection product on homology which is defined using
a count of certain pseudo-holomorphic
spheres\footnote{Pseudo-holomorphic spheres are $(j,J)$-holomorphic maps $(\C P^1,j)\to (M,J)$ for the standard complex structure $j$ on $\C P^1$ and an almost complex structure $J$ on $M$ compatible with the symplectic form.} in $M$. In fact, there are several possible algebraic setups for
this structure -- we refer the reader to
\cite{McD-Sal-psh-book}, as well as \cite{EP-toric},
\cite{Usher-deformed}, for the
precise definitions and more details.
In any case, the resulting
algebraic object is a ring with unity given by the fundamental class
$[M]$. Passing, if needed (depending on $M$ and the algebraic setup
of the construction), to an appropriate subring with unity one gets
a finite-dimensional commutative algebra with unity over a certain field that
we will denote by $\cK$. Typically, $\cK$ is the field of semi-infinite (Laurent-type) power series with coefficients in a base field $\cF$, where $\cF$ is one of the fields $\Z_p$, $\Q$, $\R$, $\C$. Abusing the terminology we will denote the
latter finite-dimensional commutative algebra by $QH (M)$ and still call it the quantum
homology of $M$.

Let us also mention the constructions of Usher \cite{Usher-deformed}
and Fukaya-Oh-Ohta-Ono \cite{FOOO-spectral-with-bulk} (the so-called {\it deformed quantum homology})
that, roughly speaking, use certain homology classes of $M$ for an additional deformation of the quantum homology product
and sometimes allow to obtain {\it different} finite-dimensional commutative algebras
as above for a given $M$ -- abusing the terminology we will still call any of these different algebras
the quantum homology of $M$ and denote it by $QH(M)$ and emphasize the difference
between them only when needed.

Given a non-zero $a\in QH(M)$ and a Hamiltonian $F: M\times
S^1\to\R$, one can define the {\it spectral number} $c(a,F)$
\cite{Schwarz}, \cite{Oh-spectral} (see
\cite{Oh1,Oh2,Viterbo-spectral} for earlier versions of the
construction and \cite{Usher-spectral,Usher-duality} for additional
important properties of the spectral numbers). It generalizes the
following classical minimax quantity: given a singular non-zero (rational)
homology class $a$ of $M$ and a continuous function $F$ on $M$,
consider the smallest value $c$ of $F$ so that $a$ can be realized
by a cycle lying in $\{ F\leq c\}$ -- for a smooth Morse function
$F$ this definition can be reformulated in terms of the Morse
homology of $F$. The construction of the spectral number $c(a,F)$ is
based on the same concept, where the singular homology is replaced
by the quantum homology of $M$ and the Morse homology of $F$ is
replaced by its {\it Floer homology}. The latter can be viewed as an
infinite-dimensional version of the Morse homology for a certain
functional, associated with $F$, on a covering of the space of free
contractible loops in $M$, with the critical points of the
functional being pre-images of contractible 1-periodic orbits of the
Hamiltonian flow of $F$ under the covering (see e.g.
\cite{McD-Sal-psh-book} for a detailed introduction to the subject).

If $F,G: M\times S^1\to \R$ are normalized and $\tphi_F = \tphi_G$,
then $c(a,F)=c(a,G)$. Thus, given $\tphi\in\tHam(M)$, one can define
$c(a,\tphi):= c(a,F)$ for any normalized $F$ generating $\tphi$.

\subsection{The main theorem}
\label{subsec-main-thm} Assume $a\in QH(M)$ is an idempotent
(for instance, $a=[M]$).
{\bf Here and further on, whenever we mention an idempotent, we assume that it is non-zero.}
Define $\mu_a : \tHam (M)\to\R$ by
\[
\mu_a (\tphi):= - {\hbox{\rm vol}\, (M)}\lim_{k\to+\infty} \frac{c(a,\tphi^k)}{k},
\]
where ${\hbox{\rm vol}\, (M)} := \int_M \omega^n$, and $\zeta_a: C^\infty (M)\to\R$ by
\[
\zeta_a (F):= \lim_{k\to+\infty} \frac{c(a,kF)}{k}.
\]
One can check \cite{EP-qst} that the limits exist and
\[
\zeta_a (F) = \frac{\int_M F\omega^n - \mu_a (\tphi_F)}{{\rm vol}\, (M)}.
\]
The next theorem shows that under certain conditions on $QH (M)$ and
$a$ the function $\mu_a$ is a homogeneous quasi-morphism and
accordingly  $\zeta_a$ is a Lie
quasi-state invariant under the adjoint action of $\tHam (M)$ on
$C^\infty (M)$, since, up to a scaling factor and an addition of a
linear map invariant under the adjoint action of $\tHam (M)$, it is
the derivative of $\mu_a$ (see Section~\ref{subsec-Lie-qs}).

We will say that $QH(M)$ is {\it field-split}, if it can be
represented, in the category of $\cK$-algebras, as a direct sum of
two subalgebras at least one of which is a field. Such a field will
be called {\it a field factor} of $QH (M)$.

\begin{thm}
\label{thm-main-field-split-case}
Assume $QH(M)$ is field-split and $a$ is the unity in a field factor of $QH(M)$. Then $\mu_a$ satisfies the following properties:

\medskip
\noindent
A. {\bf (Stability)} $\int_0^1 \min_M (F_t-G_t)\;dt \leq \mu_a(\tphi_G)-\mu_a(\tphi_F)
\leq \int_0^1 \max_M(F_t-G_t)\;dt$.

\medskip
\noindent
B. The function $\mu_a: \tHam (M)\to\R$ is a homogeneous quasi-morphism, that is,

\smallskip
\hangindent=0.6cm\hangafter=0\noindent
B1. {\bf (Homogeneity)} $\mu_a (\tphi^k) = k\mu_a (\tphi)$ for any $\tphi\in\tHam(M)$ and $k\in\Z$.

\smallskip
\hangindent=0.6cm\hangafter=0\noindent B2. {\bf (Quasi-additivity)}
There exists $C>0$ such that $|\mu_a (\tphi\tpsi)-\mu_a
(\tphi)-\mu_a (\tpsi)|\leq C$ for any $\tphi,\tpsi\in \tHam (M)$.

\medskip
\hangindent=0cm\hangafter=0\noindent C. {\bf (Calabi property)} If
$U\subset M$ is stably displaceable and $\supp F\subset U$, then $\mu_a
(\tphi_F) = \int_0^1 \int_U F_t\omega^n dt$. In other words, the
Calabi homomorphism $\tCal_U$ descends from $\tHam (U)$ to
$\tG_U\subset \tHam (M)$\footnote{This was tacitly assumed in
\cite{EP-qmm}.} and $\mu|_{\tG_U}=\tCal_U$.

\medskip
At the same time, $\zeta_a$ satisfies the following properties:

\medskip
\noindent a. {\bf (Monotonicity)} $\min_M (F-G) \leq \zeta_a
(F)-\zeta_a (G) \leq \max_M (F- G)$ for any $F,G\in C^\infty (M)$
and, in particular, if $F\leq G$, then $\zeta_a (F)\leq \zeta_a
(G)$.  Hence, $\zeta_a$ is 1-Lipschitz with respect to the uniform
norm and extends to a functional on $C(M)$ that we will still denote
by $\zeta_a$.

\medskip
\noindent
b. The functional $\zeta_a$ is a Lie quasi-state, that is

\smallskip
\hangindent=0.6cm\hangafter=0\noindent
b1. {\bf (Homogeneity)} $\zeta_a (\lambda F) = \lambda\zeta_a(F)$ for any $F\in C(M)$ and $\lambda\in\R$.

\smallskip
\hangindent=0.6cm\hangafter=0\noindent b2. {\bf (Strong
quasi-additivity)} If $F,G\in C^\infty (M)$ and $\{ F,G\} =0$, then
$\zeta_a (F+ G)= \zeta_a (F) + \zeta_a (G)$. In fact, $\zeta_a$
satisfies a stronger property: for any $F,G\in C^\infty (M)$ one has
\[
|\zeta_a (F+G)-\zeta_a (F)-\zeta_a (G)|\leq \sqrt{2C||\{F,G\}||},
\]
where $C>0$ is the constant from B2.

\medskip
\hangindent=0cm\hangafter=0\noindent c. {\bf (Vanishing property)}
If $U\subset M$ is stably displaceable and $F\in C(M)$ with $\supp F\subset
U$, then $\zeta_a (F)=0$.

\medskip
\noindent
d. {\bf (Normalization)} $\zeta_a(1)=1$.

\medskip
\noindent e. {\bf (Invariance)} $\zeta_a : C(M)\to\R$ is invariant
under the action of $Symp_0 (M)$ on $C(M)$.

\end{thm}

The functionals $\mu_a: \tHam (M)\to\R$ and $\zeta_a: C(M)\to\R$
satisfying the properties listed in
Theorem~\ref{thm-main-field-split-case} are called, respectively, a
{\it Calabi quasi-morphism} and a {\it symplectic quasi-state}. In
particular, the restriction of a symplectic quasi-state to $C^\infty
(M)$ is always a Lie quasi-state. Moreover, one can readily check
that any symplectic quasi-state is also a topological quasi-state.
The converse is true only if ${\rm dim}\, M=2$ \cite{EP-qst}.
The quasi-measure associated to a symplectic quasi-state is $\Symp_0 (M)$-invariant and vanishes
on stably displaceable open sets.

For an arbitrary $M$ and an arbitrary idempotent $a\in
QH(M)$ (for instance, $a=[M]$) one gets a weaker set of properties of $\mu_a$ and $\zeta_a$.

\begin{thm}
\label{thm-main-gen-case} Assume $a\in QH(M)$ is an arbitrary
 idempotent. Then $\mu_a$ satisfies the properties (A) and
(C) from Theorem~\ref{thm-main-field-split-case} and a weaker
version of the properties (B1) and (B2):

\medskip
\noindent B1'. {\bf (Partial homogeneity)} $\mu_a (\tphi^k) = k\mu_a (\tphi)$ for any
$\tphi\in\tHam(M)$ and $k\in\Z_{\geq 0}$.

\medskip
\noindent B2'. {\bf (Partial quasi-additivity)} Given a displaceable open
set $U\subset M$, there exists $C = C(\mu_a, U)>0$, so that $|\mu_a (\tphi\tpsi)-\mu_a
(\tphi)-\mu_a (\tpsi)|\leq C\min\{\|\tphi\|_U, \|\tpsi\|_U\}$
for any $\tphi,\tpsi\in \tHam (M)$.

\medskip
At the same time, $\zeta_a$ satisfies the properties
(a),(c),(d),(e) from Theorem~\ref{thm-main-field-split-case} and a
weaker version of the properties (b1) and (b2):

\medskip
\noindent b1'. {\bf (Partial homogeneity)} $\zeta_a (\lambda F) = \lambda\zeta_a(F)$ for any
$F\in C(M)$ and $\lambda\in\R_{\geq 0}$.

\medskip
\noindent b2'. {\bf (Partial strong quasi-additivity)} If $F,G\in
C^\infty (M)$ and $\{ F,G\} =0$ and either $\supp G$ is displaceable
or $G$ is constant, then $\zeta_a (F+ G)= \zeta_a (F) + \zeta_a
(G)$. In fact, $\zeta_a$
satisfies a stronger property: for any $F,G\in C^\infty (M)$ one has
\[
|\zeta_a (F+G)-\zeta_a (F)-\zeta_a (G)|\leq \sqrt{2C||\{F,G\}||},
\]
where $C\in (0,+\infty]$ is a constant depending on $\zeta_a$ and on the supports of $F$ and $G$.
In particular, the constant $C$ is finite if at least one of the supports is displaceable.

\end{thm}

The functionals $\mu_a$ and $\zeta_a$ satisfying the properties
listed in Theorem~\ref{thm-main-gen-case} are called, respectively,
a {\it partial Calabi quasi-morphism} and a {\it partial symplectic
quasi-state}. Clearly, a genuine Calabi quasi-morphism or a genuine
symplectic quasi-state is also a partial one.

Theorems~\ref{thm-main-field-split-case} and \ref{thm-main-gen-case}
were proved under various additional restrictions on $M$ in
\cite{EP-qmm}, \cite{EP-qst}. In
\cite{Ostr-Calabi,McDuff-monodromy,Usher-spectral,Usher-duality}
(see also \cite{EP-toric, EP-rigid}) the restrictions were removed
and the theorems were proved for new classes of closed symplectic
manifolds. The stronger part of the property (b2) in
Theorems~\ref{thm-main-field-split-case} was proved in \cite{EPZ}.
The stronger part of the property (b2') in
Theorems~\ref{thm-main-gen-case} was proved in \cite{P-quantum-noise}.
The Calabi and vanishing properties (C) and (c) in Theorems~\ref{thm-main-field-split-case}, \ref{thm-main-gen-case}
were originally proved for a displaceable $U$ -- it was later observed by Borman
\cite{Borman-reduction} that the displaceability assumption on $U$  can be weakened to stable displaceability.

Examples of closed manifolds with a field-split $QH(M)$ (for an
appropriate algebraic setup of $QH(M)$) include complex projective
spaces (or, more generally, complex Grassmanians and symplectic
toric manifolds), as well as blow-ups of symplectic manifolds -- all
with appropriate (and, in a certain sense, generic \cite{Usher-deformed})
symplectic structures
\cite{EP-qmm,EP-toric,FOOO-spectral-with-bulk,Galkin,Maydanskiy-Mirabelli,McDuff-monodromy,Ostr-Calabi,Ostrover-Tyomkin,Usher-deformed}.
The direct product of symplectic manifolds with field-split quantum
homology algebras also has this property -- possibly under some
additional assumptions on the manifolds, depending on the algebraic
setup of the quantum homology \cite{EP-toric}.
As an example of $M$ whose quantum homology $QH(M)$ is {\it not} field-split
one can take any symplectically aspherical manifold -- in such a case there are no pseudo-holomorphic
spheres in $M$ and hence there is
no difference between quantum and singular homology.
Let us emphasize
that, in general, the question whether $QH(M)$ is field-split
may depend not only on $M$ but also on the algebraic setup of the
quantum homology (and on a choice of the deformation in case of the deformed quantum homology).

\begin{exam}[\cite{EP-qmm}]
\label{exam-median-qstate} Assume $M=S^2$ and $a=[S^2]$. Then $\zeta_a: C(S^2)\to\R$ is a symplectic quasi-state and
its restriction to the set of smooth Morse functions on $S^2$ (which
is dense in $C(S^2)$ in the uniform norm) can be described in
combinatorial terms.

Namely, assume that the symplectic (that is, area) form $\omega$ on
$S^2$ is normalized so that the area of $S^2$ is $1$ and let $F$ be
a smooth Morse function on $S^2$. Consider the space $\Delta$ of
connected components of the level sets of $F$ as a quotient space of
$S^2$.
As a topological space $\Delta$ is homeomorphic to a tree. The
function $F$ descends to $\Delta$ and the push-forward of the
measure defined by $\omega$ on $S^2$ yields a non-atomic Borel
probability measure on $\Delta$. There exists a unique point
$x\in\Delta$ such that each connected component of $\Delta\setminus
x$ has measure $\leq 1/2$ (such a point $x$ is called {\it the
median} of the measured tree $\Delta$). Then $\zeta_a (F) = F(x)$.

\end{exam}

Let us note that the symplectic quasi-state $\zeta_a: C(S^2)\to\R$
in Example~\ref{exam-median-qstate} is {\it dispersion-free}, that
is, satisfies $\zeta_a (F^2) = (\zeta_a (F))^2$. Equivalently, the
corresponding quasi-measure takes only values $0$ and $1$. The
following open question is of utmost importance for the study of
symplectic quasi-states:

\begin{question}
\label{quest-dispers-free} Is it true that the 
symplectic
quasi-states constructed in Theorems~\ref{thm-main-field-split-case}
and \ref{thm-main-gen-case} are always dispersion-free?
\end{question}

\begin{rem}
\label{rem-descent} Sometimes the (partial) Calabi
quasi-morphism $\mu_a: \tHam (M)\to\R$ descends to $\Ham (M)$ (that
is, vanishes on $\pi_1 \Ham (M)$). Abusing the notation we will
denote the resulting (partial) Calabi quasi-morphism on $\Ham (M)$
also by $\mu_a$. The list of manifolds for which $\mu_a$ is known to
descend to $\Ham (M)$ for {\it all} $a$ includes symplectically aspherical manifolds
\cite{Schwarz}, complex projective spaces \cite{EP-qmm} and their
monotone products \cite{EP-qmm,Branson}, a monotone blow-up of $\C
P^2$ at three points and the complex Grassmannian $Gr (2,4)$
\cite{Branson}. The list of manifolds for which it is known that
$\mu_a$ does {\it not} descend to $\Ham (M)$ at least for {\it some} $a$ includes various
symplectic toric manifolds and,
in particular, the monotone blow-ups of $\C P^2$ at one or two
points \cite{EP-rigid,Ostr-Calabi}.

Let us note that if $(M,\omega)$ is monotone, the restriction of
$\mu_a$ to $\pi_1 \Ham (M)$ does not depend on the choice of $a$ (for a fixed algebraic setup of
$QH (M)$)
\cite{EP-rigid}. This is not necessarily true if $M$ is not monotone
\cite{Ostrover-Tyomkin}.

\end{rem}

\begin{rem}
\label{rem-non-linearity} For a closed connected $M$ the group $\Ham
(M)$ is simple and the group $\tHam (M)$ is perfect \cite{Banyaga}.
Hence, these groups do not admit non-trivial homomorphisms to $\R$
and therefore partial Calabi quasi-morphisms on $\tHam (M)$ and
$\Ham (M)$ are never homomorphisms (they are non-trivial because of
the Calabi property). Also, partial symplectic quasi-states are
never linear (use a partition of unity with displaceable supports to
check it).
Moreover, in certain cases one can verify that a partial symplectic
quasi-state $\zeta_a$ is not a genuine quasi-state (and,
accordingly, $\mu_a$ is not a genuine quasi-morphism) -- see
Section~\ref{subsec-rigidity-of-intersec}.

\end{rem}

\begin{rem}
\label{rem-Calabi-patch-up} Denote by $\cE$ the collection of all
open displaceable $U\subset (M,\omega)$ such that $\omega|_U$ is
exact. For any $U\in\cE$ the Calabi homomorphism  $\Cal_U : G_U\to
\R$ is well-defined and, by Banyaga's theorem \cite{Banyaga}, the
group $\Ker \Cal_U$ is simple, meaning that, up to a scalar factor,
$\Cal_U$ is the unique non-trivial $\R$-valued homomorphism on $G_U$
continuous on 1-parametric subgroups. If $U,V\in\cE$ and $U \subset V$,
then $G_U \subset G_V$ and
$\Cal _U = \Cal _V$ on $G_U$. Thus, if a partial Calabi
quasi-morphism $\mu_a$ descends to $\Ham (M)$, we get the following
picture: there is a family $\cE$ of subgroups of $\Ham (M)$, with
each subgroup $G_U\in\cE$ carrying the unique $\R$-valued
homomorphism $\Cal_U$ (continuous on 1-parametric subgroups), and while it is impossible to patch up all
these homomorphisms into an $\R$-valued homomorphism on $\Ham (M)$,
it is possible to patch them up into a partial Calabi
quasi-morphism $\mu_a$ (which may be non-unique).
\end{rem}

Now let us discuss  the existence and uniqueness of genuine Calabi quasi-mor\-phisms and
symplectic quasi-states on a given symplectic manifold and, in particular, the dependence of $\mu_a$ and
$\zeta_a$ on $a$ and the algebraic setup of $QH(M)$.

The set of idempotents in $QH(M)$ carries a partial order: namely,
given  idempotents $a,b\in QH(M)$, we write $a\succeq b$ if $ab=b$.
Clearly, $[M]\geq b$ for any idempotent $b\in QH(M)$. If $a\succeq
b$, then $a-b$ is also an idempotent and $a\succeq a-b$. Conversely,
if $b,b'\in QH(M)$ are two idempotents such that $bb'=0$, then
$b+b'$ is also an idempotent and $b+b'\succeq b,b'$.

The following theorem follows from basic properties of spectral numbers (cf. \cite{EP-rigid}, Theorem 1.5).

\begin{thm}
\label{thm-diff-idemp-order} Assume $a,b\in QH(M)$ are idempotents, so that $a\succeq b$. Then

a. $\mu_a\leq \mu_b$, $\zeta_a\geq \zeta_b$.

b. If $\mu_a$ is a genuine (i.e. not only partial) Calabi
quasi-morphism, then $\mu_a=\mu_b$, $\zeta_a=\zeta_b$ and thus
$\mu_b$ is also a genuine Calabi quasi-morphism and $\zeta_b$ is a
genuine symplectic quasi-state.
\end{thm}

At the same time it is possible that $a\succeq b$, $\mu_a$ is a partial Calabi quasi-morphism while $\mu_b$ is a genuine Calabi quasi-morphism -- see Examples~\ref{exam-lagr-3}, \ref{exam-exotic-torus-superheavy}.

In fact, different idempotents may define linearly independent
Calabi quasi-morphisms and symplectic quasi-states. For instance, if
$M$ is a blow-up of $\C P^2$ at one point with an appropriate
non-monotone symplectic structure, one can find Calabi quasi-morphisms $\mu_a$ and $\mu_b$ on $\tHam
(M)$ (for some idempotents $a,b\in QH(M)$) that have linearly
independent restrictions to $\pi_1\Ham(M)$ \cite{Ostrover-Tyomkin},
and if $M$ is $S^2\times S^2$ with a monotone symplectic structure,
one can find linearly independent\footnote{Note that
Calabi quasi-morphisms and symplectic quasi-states on a
given manifold form convex sets respectively in
the vector spaces of all homogeneous quasi-morphisms and all Lie quasi-states
and their linear dependence is considered
in these vector spaces.} symplectic quasi-states $\zeta_a$,
$\zeta_b$ on $C(M)$ -- see Example~\ref{exam-exotic-torus-superheavy}.

Moreover, a change of the algebraic setup of $QH(M)$ may yield new
Calabi quasi-morphisms and symplectic quasi-states on the same $M$.
For instance, if $M=\C P^n$, one can choose algebraic setups of
$QH(\C P^n)$, both for $\cF=\Z_2$ and $\cF=\C$, so that {\it in both cases}
$QH(\C P^n)$ is a field and thus the {\it only}
idempotent in $QH(\C P^n)$ is the unity $a=[\C P^n]$.
By Theorem~\ref{thm-main-field-split-case}, in both cases $a=[\C P^n]$ defines a Calabi
quasi-morphism $\mu_a$ on $\tHam (\C P^n)$ that descends to $\Ham
(\C P^n)$ (for $\cF=\C$ this is proved in \cite{EP-qmm}, the same proof works also for $\cF=\Z_2$).
However, it follows from
\cite{Wu} for $n=2$ and from \cite{Oakley-Usher} for $n=3$  that the symplectic quasi-states
defined by $[\C P^n]$ in both cases are {\it different} (see
Example~\ref{exam-RPn-Wu}).

\begin{rem}
\label{rem-cpn-givental}
Let us note that historically the first Calabi quasi-morphism on
$\tHam (\C P^n)$ was implicitly constructed by Givental in
\cite{Givental1}, \cite{Givental2} in a completely different way;
the fact that it is indeed a Calabi quasi-morphism was proved by Ben
Simon \cite{BenSimon}  (the stability property of the quasi-morphism is proved in \cite{Borman-Zapolsky}). Givental's Calabi quasi-morphism descends
from $\tHam (\C P^n)$ to $\Ham (\C P^n)$ \cite{Shelukhin}.
It would be interesting to find out whether this quasi-morphism on
$\Ham (\C P^n)$ can be expressed as $\mu_a$ for some $a\in QH(\C P^n)$
(for some algebraic setup of $QH(\C P^n)$).
\end{rem}

\begin{question}
Is the quasi-morphism $\mu_a$ for $a=[\C P^1]$ the only Calabi
quasi-morphism on $\Ham (\C P^1)$?
\end{question}

Let us note that in the case of $\C P^1$ the symplectic
quasi-state $\zeta_a$ for $a=[\C P^1]$, described in
Example~\ref{exam-median-qstate}, is known to be the unique
symplectic quasi-state on $C(S^2)$ \cite{EP-qmm}.

In some cases the non-uniqueness can be prove by using the different algebras $QH(M)$ appearing in
the deformed quantum homology construction -- see \cite{FOOO-spectral-with-bulk} for examples of
symplectic manifolds with infinitely many linearly independent
Calabi quasi-morphisms and symplectic quasi-states constructed in this way.

Let us also mention a construction due to Borman
\cite{Borman-reduction},\cite{Borman-moment-map} (also see \cite{Borman-Zapolsky}) that allows in
certain cases to use a Calabi quasi-morphism on $\tHam (N)$ in order
to build a Calabi quasi-morphism on $\tHam (M)$ if $M$ is obtained
from $N$ by a symplectic reduction or if $M$ is a symplectic
submanifold of $N$. Using different presentations of a symplectic manifold as a symplectic reduction
one can construct examples of $M$
with infinitely many linearly independent Calabi quasi-morphisms on
$\tHam (M)$ (the corresponding symplectic quasi-states are also
linearly independent).

Let us also note that apart from the Calabi quasi-morphisms from
Theorem~\ref{thm-main-field-split-case}, the Givental quasi-morphism mentioned in Remark~\ref{rem-cpn-givental} and the quasi-morphisms produced from them by Borman's reduction method, there are no known examples of
homogeneous quasi-morphisms on $\tHam (M)$ satisfying the stability property (A)
from Theorem~\ref{thm-main-field-split-case},
though otherwise there are many homogeneous
quasi-morphisms on $\tHam (M)$ (that sometimes descend to $\Ham (M)$)  -- see
\cite{E-comm-length,Gambaudo-Ghys,Py-AnnENS,Py-CR-torus,
Shelukhin-qmms}. Let us also note that for symplectic manifolds of dimension greater
than $2$
the constructions above are the only currently known constructions of partial symplectic quasi-states on closed manifolds.
(As it was
mentioned above, in dimension $2$ any topological quasi-state is
also symplectic. There is a number of ways to construct topological
quasi-states in any dimension -- see e.g.
\cite{Aar-Fundam,Knudsen-Advances}). The basic
open case for the existence of a (stable) Calabi quasi-morphism on $\tHam (M)$
is $M=T^2$ and for symplectic quasi-states it is $M=T^4$.
It is also unknown whether any
symplectic quasi-state has to come from a Calabi quasi-morphism and
whether different Calabi quasi-morphisms may define the same
symplectic quasi-state.

\begin{rem}
\label{rem-open-case}
There is a straightforward extension of the notion of a (partial)
Calabi quasi-morphism to the case of an open symplectic manifold $M$
and there are several constructions of such quasi-morphisms.

First, one can consider a conformally symplectic
embedding\footnote{A map $f:(M,\omega)\to (N,\Omega)$ between
symplectic manifolds is called {\it conformally symplectic} if
$f^*\Omega = c\omega$ for some non-zero constant $c$.} of $M$ in a
closed symplectic manifold $N$ carrying a Calabi quasi-morphism
defined on $\Ham (N)$ and use the homomorphism $\Ham (M)\to\Ham (N)$
induced by the embedding to pullback the quasi-morphism from $\Ham
(N)$ to $\Ham (M)$. (Of course, one then has to prove that the
resulting quasi-morphism on $\Ham (M)$ is non-trivial). In this way
one can, for instance, use conformally symplectic embeddings of a standard
round ball $B^{2n}$ in $\C P^n$ in order to construct
a continuum of linearly independent Calabi quasi-morphisms on $\Ham (B^{2n})$
\cite{BiEP}.

There are also {\it intrinsic} constructions of (partial) Calabi
quasi-morphisms for certain open symplectic manifolds following the
lines of the construction presented above -- see
\cite{Lanzat,Monzner-Vichery-Zapolsky} for more details as well as
for extensions of the notion of a (partial) symplectic quasi-state to the open
case. These constructions allow to extend many of the results
mentioned in this survey to the open case.
\end{rem}

\section{Applications}

\subsection{Quasi-states and rigidity of symplectic
intersections} \label{subsec-rigidity-of-intersec} A key phenomenon
in symplectic topology is rigidity of intersections of subsets of
symplectic manifolds: namely, sometimes a
subset $X$ of a symplectic manifold $M$ cannot be displaced from a
subset $Y$ by $\Ham (M)$ (or $\Symp_0 (M)$, or $\Symp (M)$), even
though $X$ can be displaced from $Y$ by a smooth isotopy. A central
role in the applications of partial symplectic quasi-states is
played by their connection to this phenomenon. Namely, to each
partial symplectic quasi-state, and, in particular, to each
idempotent $a\in QH(M)$, one can associate a certain hierarchy of
non-displaceable sets in $M$. The interplay between the hierarchies
associated to different $a$ is an interesting geometric phenomenon
in itself.

The key definitions describing the hierarchy are as follows
\cite{EP-rigid}. Let $\zeta: C(M)\to\R$ be a partial symplectic
quasi-state. We say that a closed subset $X \subset M$ is {\it heavy
with respect to $\zeta$}, if $\zeta (F) \geq \inf_X F$ for all $F\in
C(M)$, and {\it superheavy with respect to $\zeta$}, if $\zeta (F)
\leq \sup_X F$ for all $F\in C(M)$. Equivalently, $X$ is superheavy
with respect to $\zeta$, if $\zeta (F)= F(X)$ for any $F\in C(M)$
which is constant on $X$. If $\zeta=\zeta_a$ for an idempotent
$a\in QH(M)$ (and a prefixed algebraic setup of $QH(M)$),
we use the terms {\it $a$-heavy} and {\it $a$-superheavy} for
the heavy and superheavy sets with respect to $\zeta_a$. Clearly, a
closed set containing a heavy/superheavy subset is itself
heavy/superheavy. The basic properties of heavy and superheavy sets
are summarized in the following theorems.

\begin{thm}[\cite{EP-rigid}]
\label{thm-a-heavy-superheavy} Heavy and superheavy sets with
respect to a fixed partial symplectic quasi-state $\zeta$ satisfy
the following properties:

a. Every superheavy set is heavy, but, in general, not vice versa.

b. The classes of heavy and superheavy sets are $\Symp_0
(M)$-invariant.

c. Every superheavy set has to intersect every heavy set. Therefore, in view of (b),
any superheavy set cannot be displaced from any heavy
set by $\Symp_0 (M)$. In particular, any superheavy set is non-displaceable by
$\Symp_0 (M)$. On the other hand, two heavy sets may be disjoint.

d. Every heavy subset is stably non-displaceable. However, it may be
displaceable by $\Symp_0 (M)$.

e. If $\zeta$ is a genuine (that is, not partial) symplectic
quasi-state, then the classes of heavy and superheavy sets are
identical: they coincide with the class of closed sets of full
quasi-measure (that is, of quasi-measure $1$) for the quasi-measure
on $M$ associated with $\zeta$.

\end{thm}

\begin{thm}[\cite{EP-rigid}]
\label{thm-superheavy-product}
Assume that $X_i$ is an $a_i$-heavy (resp. $a_i$-superheavy)
subset of a closed connected symplectic manifold $M_i$ for some
idempotent $a_i\in QH(M_i)$, $i=1,2$. Then the product $X_1 \times
X_2\subset M_1\times M_2$ is $a_1\otimes a_2$-heavy (resp.
$a_1\otimes a_2$-superheavy)\footnote{There is an analogue of
 the K\"unneth formula for quantum homology --
in particular, to each pair of idempotents $a_1\in QH(M_1)$, $a_2\in
QH(M_2)$, one can associate an idempotent $a_1\otimes a_2\in QH
(M_1\times M_2)$. Let us note that even if $\zeta_{a_1}$ and
$\zeta_{a_2}$ are genuine symplectic quasi-states,
$\zeta_{a_1\otimes a_2}$ may be only a partial one -- see
Example~\ref{exam-exotic-torus-superheavy}.}.
\end{thm}

Changing an idempotent $a$ or changing the algebraic setup of $QH(M)$
may completely change the
heaviness/super\-heavi\-ness property of a set: there are
examples of disjoint sets that are superheavy with respect to
different idempotents in $QH(M)$ (see
Example~\ref{exam-exotic-torus-superheavy})
and there is an example of a set that is $[M]$-superheavy, if $QH(M)$ is set up over $\cF=\Z_2$, and
is disjoint from an $[M]$-superheavy set, if $QH(M)$ is set up over $\cF=\C$ (see Example~\ref{exam-RPn-Wu}).

The partial
order on the set of idempotents mentioned in
Section~\ref{sec-Calabi-qmms-sympl-qst} yields the following relation between
the corresponding collections of heavy and superheavy sets which follows
immediately from Theorem~\ref{thm-diff-idemp-order} and the examples below.

\begin{thm}
\label{thm-heavy-superheavy-diff-idemp} Assume $a,b\in QH(M)$ are
idempotents and $a\succeq b$. Then

a. Every $a$-superheavy set is also $b$-superheavy (but not necessarily vice versa).

b. Every $b$-heavy set is also $a$-heavy (but not necessarily vice versa).

\end{thm}

\begin{rem}
There is a natural action of $\Symp (M)$ on $QH (M)$. The subgroup
$\Symp_0 (M)\subset \Symp (M)$ acts on $QH(M)$ trivially and this
explains the $\Symp_0 (M)$-invariance of the partial symplectic quasi-states $\zeta_a$
and, accordingly, of the classes of $a$-heavy and
$a$-superheavy sets. If an idempotent $a\in QH(M)$ is invariant under
the action of the full group $\Symp (M)$ (e.g. if $a=[M]$), then
$\Symp_0 (M)$ can be replaced by $\Symp (M)$
everywhere in Theorems~\ref{thm-main-field-split-case},
\ref{thm-main-gen-case}, \ref{thm-a-heavy-superheavy}. In
particular, {\it any} $[M]$-superheavy set cannot be displaced from
{\it any} $a$-heavy set (for {\it any} idempotent $a\in QH(M)$) by
 $\Symp (M)$.
\end{rem}

\begin{rem}
The reason why every superheavy set $X$ (with respect to a partial
quasi-state $\zeta$) must intersect every heavy set $Y$ is very
simple: If $X\cap Y=\emptyset$, pick a function $F$ so that
$F|_X\equiv 0$ and $F|_Y\equiv 1$. Then, by the definition of heavy
and superheavy sets, $\zeta (F)=0$ and $\zeta (F)\geq 1$, which is
impossible.
\end{rem}

\begin{rem}
\label{rem-Lagr-Floer-vs-quasi-states}
For certain Lagrangian submanifolds $X$ and $Y$ of $M$ one can prove
the non-displaceability of $X$ from $Y$ by means of the Lagrangian
Floer homology of the pair $(X,Y)$ (see e.g. \cite{FOOO-book}).
The advantage of this method is that, unlike Theorem~\ref{thm-a-heavy-superheavy},
it gives a non-trivial lower bound on the number of transverse
intersection points of $\phi(X)$, $\phi\in\Ham (M)$, and $Y$. On the other hand, unlike the
Lagrangian Floer theory, symplectic quasi-states allow to prove non-displaceability results
for singular sets (see the examples below).
\end{rem}

The most basic examples of heavy and superheavy sets are an equator
in $S^2$ (that is, an embedded circle dividing $S^2$ into two parts
of equal area) which is $[S^2]$-superheavy and a meridian in $T^2$
which is $[T^2]$-heavy but not $[T^2]$-superheavy, since it is
displaceable by $\Symp_0 (T^2)$ -- in particular, it implies that
the partial symplectic quasi-state on $C(T^2)$ defined by $[T^2]$ is
not a genuine symplectic quasi-state \cite{EP-rigid} (cf.
Remark~\ref{rem-non-linearity}). The union of a meridian and a
parallel in $T^2$ is $[T^2]$-superheavy \cite{Kawasaki}. More
complicated examples come from the following constructions.

Let ${\mathbb A} \subset
C^{\infty}(M)$ be a finite-dimensional Poisson-com\-mu\-ta\-tive
vector subspace (meaning that $\{F,G\}=0$ for any $F,G\in {\mathbb A}$) . The map
$\Phi: M \to {\mathbb A}^*$ defined by  $\langle\Phi(x),F\rangle = F(x)$
is called {\it the moment map} of ${\mathbb A}$. As an example of such a moment map
one can consider the moment map of a
Hamiltonian torus action on $M$, or a map $M\to\R^N$ whose components have
disjoint supports.

A non-empty fiber $\Phi^{-1}(p)$ is called a {\it
stem} of ${\mathbb A}$ (see \cite{EP-qst}), if all non-empty fibers
$\Phi^{-1}(q)$ with $q \neq p$ are displaceable, and a {\it stable
stem}, if they are stably displaceable. If a subset of $M$ is a
(stable) stem of {\it some} finite-dimensional Poisson-commutative
subspace of $C^\infty (M)$, it will be called just {\it a (stable)
stem}. Any stem is a stable stem but possibly not vice versa
\footnote{There are no known examples of a stable stem that is not a stem, i.e. not a stem of {\it any} ${\mathbb A}$.}.

\begin{thm}[\cite{EP-qst,EP-rigid}]
\label{thm-stable-stem-superheavy}
A stable stem is superheavy with respect to any partial symplectic
quasi-state $\zeta$ on $C(M)$.
\end{thm}

Using the partial symplectic quasi-state $\zeta_{[M]}$ we get

\begin{cor}[\cite{EP-qst}]
\label{cor-existence-of-non-displ-fiber-moment-map}
For any finite-dimensional Poisson-commutative subspace of $C^{\infty}(M)$ its moment map $\Phi$ has at least one
non-displaceable fiber.
\end{cor}

The following question is closely related to Question~\ref{quest-dispers-free}.

\begin{question}
Is it true that for any finite-dimensional Poisson-commutative
subspace of $C^{\infty}(M)$ its moment map $\Phi$
has at least one {\it heavy} fiber (at least with respect to {\it some} symplectic quasi-state $\zeta$ on $C(M)$)?
\end{question}

\begin{rem}
\label{rem-pushforward-of-quasi-measures}
If $\zeta$ a genuine symplectic quasi-state on $C(M)$, then Theorem~\ref{thm-stable-stem-superheavy}
can be proved using the quasi-measure $\tau$ associated to $\zeta$ \cite{EP-qst}.
Namely, the push-forward of $\tau$ to ${\mathbb A}^*$ by
the moment map $\Phi$ of a Poisson commutative subspace ${\mathbb A}$
is a Borel probability measure $\nu$ on ${\mathbb A}^*$ \cite{EP-qst}.
As we already noted above, the vanishing property of $\zeta$
implies that $\tau$ vanishes on stably displaceable open subsets of $M$.
Therefore if a fiber $\Phi^{-1} (p)$ of $\Phi$ is a stable stem, the support of $\nu$ has to
be concentrated at $p$, meaning that $\tau (\Phi^{-1} (p))=1$ or, in other words (see Theorem~\ref{thm-a-heavy-superheavy}), $\Phi^{-1} (p)$ is superheavy with respect to $\zeta$.
\end{rem}

Here are a few examples of (stable) stems \cite{EP-qst,EP-rigid}.
The Lagrangian Clifford torus in $\C P^n$, defined as $L=\{
[z_0:\ldots :z_n]\in \C P^n\ |\  |z_0|=\ldots=|z_n|\}$, is a stem,
hence $[\C P^n]$-superheavy (this generalizes the example of an
equator in $S^2$ mentioned above). The codimension-1 skeleton of a
triangulation of a closed symplectic manifold $M^{2n}$ all of whose
$2n$-dimensional simplices are displaceable is a stem. A fiber $\Phi^{-1} (0)$
of the normalized moment map $\Phi$ of a
compressible\footnote{An effective Hamiltonian $T^k$-action on
$(M,\omega)$ is called {\it compressible} if the image of the
homomorphism $\pi_1(T^k)\to \pi_1(\Ham(M))$, induced by the action,
is a finite group.} Hamiltonian torus action on $M$ is a stable stem\footnote{Stable stems appearing in this way potentially may not be
genuine stems.}. A direct product of (stable) stems is a (stable)
stem and that the image of a (stable) stem under {\it any}
symplectomorphism is again a (stable) stem.

In case when the Hamiltonian torus action is not compressible, much
less is known about (stable) displaceability of fibers of the moment
map of the action (aside from the case of symplectic toric manifolds where many
results have been obtained in recent years by different authors). If
$(M,\omega)$ is monotone, one can explicitly find a
so-called {\it special fiber} of the moment map $\Phi$ of the action
which is $a$-superheavy for {\it any} idempotent $a\in
QH(M)$ \cite{EP-rigid}. For a Hamiltonian $T^n$-action torus on a
monotone $(M^{2n},\omega)$ (that is, for a monotone symplectic toric
manifold) the special fiber (which in this case is a Lagrangian
torus) can be described in simple combinatorial terms involving the
moment polytope (that is,
the image of the moment map which is a convex polytope) -- see \cite{EP-rigid}. It is not known
whether in the latter case the special fiber is always a stem -- see
\cite{McDuff-probes} for a detailed
investigation of this question. Interestingly enough, the question
whether the special fiber of the {\it normalized} moment map for a
monotone symplectic toric manifold $M$ coincides with the fiber over
zero is related to the existence of a K\"ahler-Einstein metric on
$M$ -- see \cite{EP-rigid,Shelukhin}.

Finally, heaviness/superheaviness of Lagrangian submanifolds can be
proved using various versions of Lagrangian Floer homology. Namely,
to certain Lagrangian submanifolds $L$ of $M$ one can associate the
{\it quantum homology} (or the {\it Lagrangian Floer homology})
$QH(L)$ that comes with an {\it open-closed map} $i_L: QH(L)\to
QH(M)$
\cite{Albers,Biran-Cornea-uniruling,FOOO-book,FOOO-spectral-with-bulk}.
If $i_L (x)$ is non-zero for certain $x\in QH(L)$, then $L$ is
$[M]$-heavy\footnote{The $a$-heaviness for an arbitrary idempotent $a\in QH(M)$ can also be proved by the same method under a stronger non-vanishing assumption.}, and if $i_L(x)$ divides an idempotent $a\in QH(M)$, then
$L$ is $a$-superheavy -- this is shown in \cite{EP-rigid} in the monotone case, cf.
\cite{Biran-Cornea-uniruling,FOOO-book,FOOO-spectral-with-bulk}.

Here are a few examples where this method can be applied. Let us
emphasize that the applications to specific $L$ do depend on a proper
choice of the algebraic setup for the quantum homology in each case -- see e.g. Example~\ref{exam-RPn-Wu};
we will ignore this issue in the other examples below and refer the reader to
\cite{Biran-Cornea-uniruling,EP-rigid,FOOO-spectral-with-bulk}
for details.

\begin{exam}[\cite{EP-rigid}] \label{exam-lagr-1}
Assume that $L\subset M$ is a Lagrangian submanifold and $\pi_2 (M,L)=0$.
Then $L$ is $[M]$-heavy. Note that in this case heaviness may not be
improved to superheaviness: the meridian in $T^2$ is $[T^2]$-heavy
but not $[T^2]$-superheavy.
\end{exam}

\begin{exam}\label{exam-RPn-Wu}
The real projective space $\R P^n$, which is a Lagrangian
submanifold of $\C P^n$, is $[\C P^n]$-superheavy, {\it as long as  $QH(\C P^n)$ is set up over $\cF=\Z_2$}
\cite{Biran-Cornea-uniruling,EP-rigid}.
In particular, this implies that $\R P^n$ is not displaceable by $\Symp (\C P^n)$ from the
Clifford torus (see \cite{Alston,Tamarkin} for other proofs of this
fact).

On the other hand, $\R P^n$ may not be $[\C P^n]$-superheavy, if $QH(\C P^n)$ is set up over $\cF=\C$ --
for $n=2$ this follows from \cite{Wu} and for $n=3$ from \cite{Oakley-Usher}. This implies that the symplectic
quasi-states defined by $[\C P^n]$ for the setups of $QH(M)$ over $\cF=\Z_2$ and $\cF=\C$ are different for $n=2,3$.
\end{exam}

\begin{exam}[\cite{EP-rigid}] \label{exam-lagr-3}
Consider the torus $T^{2n}$ equipped with the standard symplectic
structure $\omega = dp\wedge dq$. Let $M^{2n}=T^{2n}\sharp
\overline{\C P^{n}}$ be a symplectic blow-up of $T^{2n}$ at one
point (the blow-up is performed in a small ball $B$ around the
point). Assume that the Lagrangian torus $L\subset T^{2n}$ given by
$q=0$ does not intersect $B$.

Then the proper transform of $L$ is a Lagrangian submanifold of $M$
which is $[M]$-heavy but {\it not} $a$-heavy for some other
idempotent $a\in QH(M)$ (that, roughly speaking, depends on the
exceptional divisor of the blow-up). In this case the functional
$\zeta_{[M]}$ is a partial (but not genuine) symplectic quasi-state
while $\zeta_a$ is a genuine symplectic quasi-state.
\end{exam}

\begin{exam}[\cite{EP-rigid, Eliashb-Polt}] \label{exam-exotic-torus-superheavy}
Let $S^2$ be the standard unit sphere in $\R^3$ with the standard
area form $\sigma$. Let $M:=S^2\times S^2$ be equipped with the
symplectic form $\sigma\oplus\sigma$. Denote by $x_1, y_1,z_1$ and
$x_2,y_2,z_2$ the Euclidean coordinates on the two $S^2$-factors.

Consider the following three Lagrangian submanifolds of $M$: the
anti-diagonal $\Delta:= \{(u,v)\in S^2 \times S^2 \ :\
u=-v\}$, the Clifford
torus $L = \{ z_1=z_2=0\}$ and the torus $K= \{ z_1+z_2=0,\  x_1 x_2 +y_1y_2 +z_1 z_2= -1/2\}$. Clearly,
$L$ intersects both $K$ and $\Delta$, while $K\cap\Delta =
\emptyset$.

For a certain algebraic setup of $QH(M)$ (with $\cF=\C$) the algebra $QH(M)$ is a
direct sum of two fields whose unities will be denoted by $a_{-}$
and $a_{+}$ (in particular, $a_{-}+a_{+}=[M]$). The idempotents
$a_{-}$ and $a_{+}$ define symplectic quasi-states $\zeta_{a_{-}}$,
$\zeta_{a_{+}}$. At the same time $\zeta_{[M]}$ is only a partial,
not genuine, symplectic quasi-state. The submanifolds $\Delta$ and $L$ are
$a_{-}$-superheavy, while $K$ is not. At the same time $K$ and $L$
are $a_+$-superheavy, while $\Delta$ is not. All three sets
$\Delta$, $K$, $L$ are $[M]$-heavy but $L$ is the only one of them
that is $[M]$-superheavy.

In particular, the Lagrangian tori $L$ and $K$ cannot be mapped into
each other by any symplectomorphism of $M$. See
\cite{
Oakley-Usher} and the
references therein for more results on the Lagrangian torus $K$.
\end{exam}

For more examples of heavy and superheavy Lagrangian submanifolds
obtained by means of an open-closed map see
\cite{EP-rigid,FOOO-spectral-with-bulk}.

\subsection{Quasi-states and connecting trajectories of Hamiltonian flows}
\label{subsec-ham-chords} Here is an application of symplectic
quasi-states to Hamiltonian dynamics. As above, we assume that $M$
is a closed connected symplectic manifold.

\begin{thm}[\cite{BEP}]\label{thm-non-autonomous-superheavy}
 Let $X_0,X_1,Y_0,Y_1 \subset M$ be a quadruple of closed sets
 so that
$X_0 \cap X_1 = Y_0 \cap Y_1 =\emptyset$ and the sets $X_0 \cup Y_0,
Y_0 \cup X_1, X_1 \cup Y_1, Y_1 \cup X_0$ are all $a$-superheavy for
an idempotent $a\in QH(M)$. Let $G \in C^\infty (M
\times S^1)$ be a $1$-periodic Hamiltonian with $G_t|_{Y_0} \leq 0$,
$G_t|_{Y_1} \geq 1$ for all $t\in S^1$.

Then there exists a point $x \in M$ and time moments $t_0,t_1 \in
\R$ so that $\phi_G^{t_0} (x) \in X_0$ and $\phi_G^{t_1} (x) \in
X_1$. Furthermore, $|t_0-t_1|$ can be bounded from above by a
constant depending only on $a$, if $G$ is
time-independent, and both on $a$ and the
oscillation $\max_{M\times S^1} G - \min_{M\times S^1} G$ of $G$,
if $G$ is time-dependent.
\end{thm}

\begin{rem}
\label{rem-pb4}
The proof of Theorem~\ref{thm-non-autonomous-superheavy} uses the following important notion \cite{BEP}: Given a quadruple
$X_0, Y_0, X_1, Y_1$ of compact subsets of a (possibly open)
symplectic manifold such that $X_0\cap X_1=Y_0\cap Y_1=\emptyset$, define
$pb_4 (X_0, Y_0, X_1, Y_1)$ (where $pb$ stands for the ``Poisson brackets" and $4$ for the number of sets) by
$pb_4 (X_0, Y_0, X_1, Y_1) = \inf_{F,G} ||\{ F,G\}||$,
where the infimum is taken over all compactly supported
smooth $F,G$ satisfying $F|_{X_0}\leq 0$, $F|_{X_1}\geq 1$, $G|_{Y_0}\leq 0$, $G|_{Y_1}\geq 1$.

Given a quadruple $X_0, Y_0, X_1, Y_1$ as in
Theorem~\ref{thm-non-autonomous-superheavy}, one can use the strong
form of the property (b2) in Theorem~\ref{thm-main-field-split-case}
to prove the positivity of $pb_4$ for certain stabilizations
\cite{BEP} of the four sets. The existence of a connecting
trajectory of the Hamiltonian flow is then deduced from the
positivity of $pb_4$ using an averaging argument \cite{BEP}.
\end{rem}

\begin{exam}[\cite{BEP}]\label{exam-S2xS2-qstates-dynamics}
Consider an open tubular neighborhood $U$ of the zero-section $T^n$ in $T^* T^n$.
Pick $q_0, q_1\in T^n$ and consider the open cotangent disks $D_i=T^*_{q_i} T^n\cap U$, $i=0,1$.
Let $M= S^2 \times \ldots \times S^2$ be the
product of $n$ copies of $S^2$ equipped with the split symplectic
structure $\omega= \sigma \oplus\ldots\oplus \sigma$. Let $Y_1\subset M$ be the product of equators in the $S^2$ factors. It is
a Lagrangian torus.
If $\int_{S^2} \sigma$ is sufficiently large, then, by the Weinstein neighborhood theorem (see \cite{McD-Sal-intro}), $U$
can be symplectically identified with a tubular neighborhood of $Y_1$ in $M$.
Using the identification we consider $X_i:= \overline{D}_i$, $i=0,1$, $Y_1$ and $U$ as subsets of $M$. Set
$Y_0:=M\setminus U$.

If $\int_{S^2} \sigma$ is sufficiently large, the quadruple $X_0,Y_0,X_1,Y_1\subset M$ satisfies the
assumptions of Theorem~\ref{thm-non-autonomous-superheavy} \cite{BEP}.
In particular, let $G: T^* T^n \times S^1 \to \R$
be a Hamiltonian supported in $U\subset T^* T^n$
which is $\geq 1$ on $T^n \times S^1$.
Theorem~\ref{thm-non-autonomous-superheavy} implies
the existence of a trajectory
of the Hamiltonian flow of $G$ passing through
$D_0$ and $D_1$.

Switching the pair $X_0, X_1$ with the pair $Y_0, Y_1$ and applying
Theorem~\ref{thm-non-autonomous-superheavy} to the switched pairs
one can show in a similar way that if $F: T^* T^n\times S^1\to\R$ is
a compactly supported Hamiltonian such that $F|_{D_0\times S^1} \leq 0$,
$F|_{D_1\times S^1} \geq 1$, then there exists a trajectory of the
Hamiltonian flow of $F$ connecting the zero-section of $T^* T^n$
with $\partial U$. It would be interesting to find out whether this
fact can be related to the well-known Arnold diffusion phenomenon in Hamiltonian dynamics
that concerns trajectories of a similar
kind.\end{exam}

\subsection{Quasi-states and $C^0$-rigidity of Poisson brackets}\label{subsec-Poisson-brackets}
In this section we still assume that $(M,\omega)$
is a closed connected symplectic manifold.

The Poisson brackets of two smooth functions on $(M,\omega)$
depend on their first derivatives. Nevertheless,
as it was first discovered in \cite{Cardin-Viterbo}, the Poisson
bracket displays a certain rigidity with respect to the {\it
uniform} norm of the functions. This rigidity is best expressed in
terms of the {\it profile function} defined as follows \cite{BEP}.

Equip the space $\Pi := C^\infty (M)\times C^\infty (M)$ with the
product uniform metric: $d ((F,G), (H,K)) = ||F-H|| +||G-K||$. For
each $s\geq 0$ define $\Pi_s := \{ (H,K)\in \Pi\ |\ ||\{ H,K\}||\leq
s\}$. In particular, $\Pi_0$ is the set of Poisson-commuting pairs.
Given a pair $(F,G)\in\Pi$, define the {\it profile function}
$\rho_{F,G}: \R_{\geq 0}\to \R_{\geq 0}$ by $\rho_{F,G} (s) =
d((F,G),\Pi_s)$.

\begin{question}
Given a pair $(F,G)\in\Pi$, what can be said of $\rho_{F,G} (0)$? In
other words, how well can $(F,G)$ be approximated with respect to $d$ by a
Poisson-commuting pair?
\end{question}

Let us note that similar approximation questions have been
extensively studied for matrices -- see e.g. \cite{Hastings} and the
references therein. It follows from \cite{EP-Poisson-single} that the sets $\Pi_s$, $s\geq 0$,
are closed with respect to $d$ and therefore
$\rho_{F,G} (s)>0$ for $s\in [0, ||\{F,G\}||)$. Symplectic
quasi-states help to give a more precise answer in certain cases.

\begin{thm}[\cite{BEP}]
\label{thm-profile-function-at-zero} Let $\zeta: C(M)\to\R$ be a
symplectic quasi-state.

a. Assume $X,Y,Z\subset M$ are closed sets that are superheavy with
respect to $\zeta$ and satisfy $X\cap Y\cap Z=\emptyset$. Assume
$F|_X\leq 0$, $G|_Y\leq 0$, $(F+G)|_Z\geq 1$ and at least one of the
functions $F,G$ has its range in $[0,1]$. Then $\rho_{F,G} (0)=1/2$
and for some positive constant $C$, independent of $F,G$, and for
all $s \in [0; ||\{F,G\}||]$
\[
\frac{1}{2}- C\sqrt{s} \leq \rho_{F,G}(s)\leq \frac{1}{2}- \frac{s}{2||\{F,G\}||}.
\]

b.  Let $X_0,X_1,Y_0,Y_1\subset M$ be closed sets so that $X_0 \cap
X_1 = Y_0 \cap Y_1 =\emptyset$ and the sets $X_0 \cup Y_0, Y_0 \cup
X_1, X_1 \cup Y_1, Y_1 \cup X_0$ are all superheavy with respect to
$\zeta$. Assume $F,G\in C^\infty (M)$, $F|_{X_0}\leq 0$,
$F|_{X_1}\geq 1$, $G|_{Y_0}\leq 0$, $G|_{Y_1}\geq 1$ and at least
one of the functions $F,G$ has its range in $[0,1]$. Then
$\rho_{F,G} (0)=1/2$ and for some positive constant $C$, independent
of $F,G$, and for all $s \in [0; ||\{F,G\}||]$
\[
 \frac{1}{2}- Cs \leq \rho_{F,G}(s)\leq \frac{1}{2}- \frac{s}{2||\{F,G\}||}.
\]
\end{thm}
The proof of part (b) uses the fact that $pb_4 (X_0,Y_0,X_1,Y_1)>0$ (see Remark~\ref{rem-pb4})
and part (a) is based on the positivity of a similar Poisson
bracket invariant $pb_3 (X,Y,Z)$ -- see \cite{BEP} for more details, as well as for
examples where the theorem can be applied, including
an example where the lower bound in part (a) is asymptotically sharp.
For a version of Theorem~\ref{thm-profile-function-at-zero} for {\it iterated}
Poisson brackets see \cite{EPR}.

Here is another fact concerning the $C^0$-rigidity of Poisson
brackets whose proof uses partial symplectic quasi-states and
the strong version of their partial quasi-additivity
(see Theorem~\ref{thm-main-gen-case}). Let
$\cU=\{U_1,\ldots,U_N\}$ be a finite cover of $M$ by
displaceable open sets. Given a partition of unity
$\vec{F}=\{F_1,\ldots,F_N\}$ subordinated to $\cU$ (that is, $\supp
F_i\subset U_i$ for every $i$), consider the following measure of
its Poisson non-commutativity:
\[
\kappa (\vec{F}) := \inf_{x,y\in [-1,1]^N} \Bigl\| \Big\{\sum_{i=1}^N x_i
F_i,\sum_{j=1}^N y_j F_j \Big\}\Bigr\|,
\]
where the infimum is taken over all $x=(x_1,\ldots,x_N),
y=(y_1,\ldots,y_N)\in [-1,1]^N$. Set $pb (\cU) :=\inf_{\vec{F}}
\kappa(\vec{F})$, where the infimum is taken over all partitions of
unity $\vec{F}$ subordinated to $\cU$. We say that $\cU$ is {\it dominated by an open set $U\subset M$} if for
each $i=1,\ldots,N$ there exists $\phi_i\in\Ham (M)$ so that
$U_i\subset \phi_i (U)$.

\begin{thm}[\cite{P-quantum-noise}]\label{thm-cover}
Assume $\cU$ is dominated by a displaceable open set $U$. Then there
exists a constant $C=C(U)>0$ so that
\begin{equation}
\label{eqn-pb-cover}
pb (\cU)\geq C/N^2.
\end{equation}
\end{thm}

The theorem strengthens a similar result proved previously in
\cite{EPZ}. It is not clear whether the inequality
\eqref{eqn-pb-cover} can be improved -- see \cite{P-quantum-noise}
for a discussion.

\begin{rem}
\label{rem-from-symplectic-to-quantum-mechanics}
Amazingly, Theorem~\ref{thm-cover} that belongs to the mathematical
formalism of classical mechanics can be used to prove results about
mathematical objects of quantum nature appearing in the
Berezin-Toeplitz quantization of a symplectic manifold -- see
\cite{P-quantum-unsharpness,P-quantum-noise,PR-book}.
\end{rem}

\subsection{Quasi-morphisms and metric properties of $\Ham (M)$}
\label{subsec-qmms-and-props-of-Ham} The group $\Ham (M)$ carries
various interesting metrics.
Here we will discuss how Calabi
quasi-morphisms can be used to study these metrics.

The most remarkable metric on $\Ham (M)$ is the Hofer metric. Namely, define
the {\it Hofer norm} $||\phi||_H$ of $\phi\in\Ham (M)$
as $||\phi||_H = \inf_F \int_0^1 (\max_M F_t - \min_M
F_t) dt$, where the infimum is taken over all (time-dependent)
Hamiltonians $F$ generating $\phi$ (if $M$ is open, $F$ is also
required to be compactly supported). The {\it Hofer metric} is
defined by $\varrho (\phi,\psi) = ||\phi\psi^{-1}||_H$. It is a deep
result of symplectic topology that $\varrho$ is a bi-invariant
metric -- see e.g. \cite{McD-Sal-intro} and the references therein.

Assume $\Ham (M)$ admits a partial Calabi quasi-morphism $\mu$. Then the
stability property of $\mu$ (see
Theorem~\ref{thm-main-field-split-case}) implies that
$\mu$ is Lipschitz with respect to $\varrho$
and therefore the diameter of $\Ham
(M)$ with respect to $\varrho$ is infinite \cite{EP-qmm}. Moreover,
the Lipschitz property of
$\mu$ with respect to $\varrho$ allows to obtain the following
result on the growth of 1-parametric subgroups of
$\Ham (M)$ with respect to the Hofer norm.

\begin{thm}[\cite{PR-book}, cf. \cite{EP-qmm}]
Assume $\Ham (M)$ admits a partial Calabi quasi-morphism. Then there exists
a set $\Xi\subset C^\infty (M)$ which is $C^0$-open and
$C^\infty$-dense in $C^\infty (M)$ so that $\displaystyle \lim_{t\to +\infty}
\frac{||\phi_F^t||_H}{t || F ||}>0$ for any
$F\in\Xi$.
\end{thm}

Calabi quasi-morphisms can be also applied to the study of the metric induced by
$\varrho$ on certain spaces of Lagrangian submanifolds of $M$ -- see \cite{Khanevsky, Seyfaddini-LagrHofer}.

The group $\Ham (M)$ also carries the $C^0$-topology: equip $M$ with
a distance function $d$, given by a Riemannian metric on $M$, assume $d$ is bounded, and
define the $C^0$-topology on $\Ham (M)$ as the one induced by the
metric
${\it dist} (\phi,\psi) = \max_{x\in M} d(\phi (x),\psi (x))$. The relation between the $C^0$-topology and
the Hofer metric on $\Ham (M)$ is rather delicate (for instance, the
$C^0$-metric is never continuous with respect to the Hofer metric).
One can use the Calabi quasi-morphisms on $\Ham (B^{2n})$ (see
Remark~\ref{rem-open-case}) in order to construct infinitely many
linearly independent homogeneous quasi-morphisms on $\Ham (B^{2n})$
that are both Lipschitz with respect to the Hofer metric and
continuous in the $C^0$-topology \cite{EPP}. This yields the
following corollary answering a question of Le Roux \cite{LeRoux}:

\begin{cor}[\cite{EPP}]
For any $c\in\R$ the set $\{ \phi\in\Ham (B^{2n})\ |\ ||\phi||_H\geq
c\}$ has a non-empty interior in the $C^0$-topology.
\end{cor}

See \cite{Seyfaddini} for an extension of this result to a wider class of open symplectic manifolds.

(Partial) Calabi quasi-morphisms can be also applied to
the study of the norms $\|\phi\|_U$ and $\|\phi\|_{U,0}$ on $\Ham (M)$
(see Section~\ref{subsec-sympl-prelim}) and the metrics defined by them.
Namely, let $U\subset M$ be a displaceable open set. Assume $\Ham (M)$
admits a Calabi quasi-morphism $\mu$. A standard fact about
quasi-morphisms bounded on a generating set yields that there exists
a constant $C=C(\mu)>0$ so that $\|\phi\|_{U,0}\geq C|\mu (\phi)|$
for any $\phi\in\Ham (M)$. In particular, $\Ham (M)$ is unbounded
with respect to the norm $\|\cdot\|_{U,0}$.

On the other hand, if $\mu: \Ham (M)\to\R$ is only a partial {\it
but not genuine} Calabi quasi-morphism, it can be used to show the
unboundedness of $\Ham (M)$ with respect to the norm $\|\cdot\|_U$
\cite{Burago-Ivanov-Polt}. Namely,  assume $X,Y=\varphi(X)\subset
M$, $\varphi\in \Symp_0 (M)$, are heavy\footnote{With respect to the
partial symplectic quasi-state $\zeta$ associated to $\mu$. If $X$
is heavy with respect to $\zeta$, then so is $\varphi(X)$, since
$\zeta$ is $\Symp_0 (M)$-invariant.} disjoint\footnote{This is
possible only if $\mu$ is {\it not} a genuine Calabi quasi-morphism,
since otherwise each heavy set (with respect to $\zeta$) is also
superheavy and therefore must intersect any other heavy set.} closed
subsets of $M$ and $V,W$ are their disjoint open neighborhoods (for
instance, $X$ and $Y$ can be two meridians on a standard symplectic
torus). Let $F$, $G$ be smooth functions supported, respectively, in
$V,W$ so that $F|_X\equiv 1= \max_M F$, $G|_Y\equiv 1= \max_M G$.
One can easily show that
\[
k=|\mu (\phi^k_F\phi^k_G) - \mu (\phi^k_F) - \mu (\phi^k_G)| \leq
C\|\phi^k_F\|_U
\]
for any $k\in\N$ and a constant $C>0$ depending only on $\mu$ and $U$. Thus in such a case
$\|\phi^k_F\|_U$ grows asymptotically linearly with $k$
and the norm $\|\cdot\|_U$ is unbounded on $\Ham (M)$.

\begin{question}
Does there exist a closed symplectic manifold $M$ for which $\|\cdot\|_U$ is bounded on $\Ham (M)$?
\end{question}

\subsection{First steps of symplectic function theory -- discussion}\label{subsec-sympl-fcn-theory}
A smooth manifold $M$ and various geometric structures on it can be described in terms of the function space $C^\infty (M)$ -- for instance, subsets of $M$ correspond to ideals in $C^\infty (M)$, tangent vectors to derivations on $C^\infty (M)$ etc. In particular, a symplectic structure on $M$ is completely determined by the corresponding Poisson brackets on $C^\infty (M)$ which means that, in principle, any symplectic phenomenon has a counterpart in the {\it symplectic function theory}, that is, the function theory of the Poisson brackets. The key feature of symplectic topology is $C^0$-rigidity appearing in various forms for smooth objects on symplectic manifolds. Its counterpart in the symplectic function theory is the rigidity of the Poisson brackets with respect to the $C^0$-norm on functions -- see e.g. \cite{PR-book}
for a deduction of the foundational Eliashberg-Gromov theorem on the $C^0$-closedness of $Symp (M)$ from the $C^0$-rigidity of the Poisson brackets.

The results and methods presented in this survey show that thinking about symplectic phenomena in terms of the function theory may have a number of advantages. First, it allows to use the ``Lie group - Lie algebra" connection between $\tHam (M)$ and $C^\infty (M)/\R$: most of the properties of symplectic quasi-states are proved using the properties of Calabi quasi-morphisms. Second, it allows to deal with singular sets (see Remark~\ref{rem-Lagr-Floer-vs-quasi-states}). Third, it allows to apply functional methods, like averaging, to Hamiltonian dynamics (see Remark~\ref{rem-pb4}). Fourth, it helps to find connections between symplectic topology and quantum mechanics since it is the Poisson algebra $C^\infty (M)$ that is being quantized in various
quantization constructions (see e.g.
Remark~\ref{rem-from-symplectic-to-quantum-mechanics}). Moreover, one may hope to discover new geometric and dynamical phenomena by studying the function theory of the Poisson bracket. For instance,
the behavior of the profile function $\rho (t)$ as $t\to 0$ (see Section~\ref{subsec-Poisson-brackets} and \cite{BEP}) is obviously of interest in symplectic function theory but its geometric or dynamical implications are absolutely unclear.

\bibliographystyle{alpha}

\begin{thebibliography}{7}

\bibitem{Aar91} Aarnes, J.F.,
Quasi-states and quasi-measures, \emph{Adv. Math.} \textbf{86}
(1991), 41-67.


\bibitem{Aar-Fundam}  Aarnes, J.F.,
Construction of non-subadditive measures and discretization of Borel measures,\emph{ Fund. Math.} \textbf{147} (1995), 213--237.

%


\bibitem{Albers} Albers, P.,
 On the extrinsic topology of Lagrangian submanifolds, \emph{Int.
Math. Res. Not.} \textbf{38} (2005), 2341-2371. Erratum: \emph{Int.
Math. Res. Not.} (2010), 1363--1369.


\bibitem{Alston} Alston, G.,
Lagrangian Floer homology of the Clifford torus and real projective space in odd dimensions, \emph{J. Sympl. Geom.}
\textbf{9} (2011), 83--106.



\bibitem{Banyaga} Banyaga, A.,
Sur la structure du groupe des diff{\'e}omorphismes qui
pr{\'e}servent une forme symplectique,  \emph{Comm. Math. Helv.} \textbf{53} (1978), 174--227.



\bibitem{Bell} Bell, J.S.,  On the problem of hidden variables in quantum mechanics,
\emph{Rev. Modern Phys.} \textbf{38} (1966), 447--452.

\bibitem{BenSimon}  Ben Simon, G., The nonlinear Maslov index and the Calabi homomorphism,
\emph{Comm. Contemp. Math.} \textbf{9} (2007), 769--780.


\bibitem{Biran-Cornea-uniruling}  Biran, P., Cornea, O., Rigidity and uniruling for Lagrangian submanifolds,
\emph{Geom. Topol.} \textbf{13} (2009), 2881--2989.

\bibitem{BiEP} Biran, P., Entov, M., Polterovich, L.,  Calabi
quasimorphisms for the symplectic ball,  \emph{Comm. Contemp. Math.}
\textbf{6} (2004),  793--802.

\bibitem{Borman-reduction} Borman, M.S., Symplectic reduction of quasi-morphisms and quasi-states, \emph{J. Sympl. Geom.}
\textbf{10} (2012), 225--246.

\bibitem{Borman-moment-map} Borman, M.S., Quasi-states, quasi-morphisms, and the moment map,
\emph{Int. Math. Res. Not.} (2013), 2497--2533.

\bibitem{Borman-Zapolsky} Borman, M.S., Zapolsky, F., Quasi-morphisms on contactomorphism groups and contact rigidity, preprint, arXiv:1308.3224, 2013.

\bibitem{Branson} Branson, M., Symplectic manifolds with vanishing action-Maslov homomorphism, \emph{Algebr. Geom. Topol.}
\textbf{11} (2011), 1077--1096.

\bibitem{BEP} Buhovsky, L., Entov, M., Polterovich, L., Poisson brackets and symplectic invariants,
\emph{Selecta Math. (N.S.)} \textbf{18} (2012), 89--157.

\bibitem{Burago-Ivanov-Polt}  Burago, D., Ivanov, S., Polterovich, L.,
Conjugation-invariant norms on groups of geometric origin, in
\emph{Groups of Diffeomorphisms: In Honor of Shigeyuki Morita on
the Occasion of His 60th Birthday}. Adv. Studies in Pure
Math.~52, Math. Soc. of Japan, Tokyo, 2008.

\bibitem{Calegary-scl}  Calegari, D., \emph{scl}. MSJ Memoirs~20, Math. Soc. of Japan, Tokyo, 2009.

\bibitem{Cardin-Viterbo} Cardin, F., Viterbo, C.,  Commuting
Hamiltonians and Hamilton-Jacobi multi-time equations, \emph{Duke Math.
J.} \textbf{144} (2008), 235--284.


\bibitem{Eliashb-Polt} Eliashberg, Y., Polterovich, L., Symplectic quasi-states on the quadric surface
and Lagrangian submanifolds, preprint, arXiv:1006.2501, 2010.

\bibitem{E-comm-length} Entov, M.,  Commutator length of
symplectomorphisms, \emph{Comm. Math. Helv.} \textbf{79} (2004), 58--104.

\bibitem{EP-qmm} Entov, M., Polterovich, L.,  Calabi quasimorphism
and quantum homology, \emph{Int. Math. Res. Not.} \textbf{30} (2003),
1635-1676.

\bibitem{EP-qst} Entov, M., Polterovich, L.,  Quasi-states and
symplectic intersections, \emph{Comm. Math. Helv.} \textbf{81} (2006), 75--99.

\bibitem{EP-toric} Entov, M., Polterovich, L., Symplectic
quasi-states and semi-simplicity of quantum homology, in \emph{Toric
Topology (eds. M.Harada, Y.Karshon, M.Masuda and T.Panov), 47--70}.
Contemp. Math.~460, AMS, Providence RI, 2008.

\bibitem{EP-rigid} Entov, M., Polterovich L.,  Rigid subsets
of symplectic manifolds, \emph{Compositio Math.} \textbf{145} (2009),
773--826.

\bibitem{EP-Lie-qs} Entov, M., Polterovich, L., Lie quasi-states,
\emph{J. Lie Theory} \textbf{19} (2009), 613--637.

\bibitem{EP-Poisson-single} Entov, M., Polterovich, L.,
$C^0$-rigidity of Poisson brackets, in \emph{Proceedings of the
Joint Summer Research Conference on Symplectic Topology and
Measure-Preserving Dynamical Systems (eds. A. Fathi, Y.-G. Oh and C.
Viterbo), 25--32}. Contemp. Math.~512, AMS,
Providence RI, 2010.

\bibitem{EPP} Entov, M., Polterovich, L., Py, P., On continuity of quasimorphisms for symplectic maps. With an appendix by Michael Khanevsky, in \emph{Perspectives in analysis, geometry, and topology, 169--197}. Progr. Math.~296, Birkhäuser/Springer, New York, 2012.

\bibitem{EPR} Entov, M., Polterovich L., Rosen, D.,  Poisson
brackets, quasi-states and symplectic integrators,  \emph{Discr. and
Cont. Dyn. Systems} \textbf{28} (2010), 1455--1468.

\bibitem{EPZ} Entov, M., Polterovich, L., Zapolsky, F.,  Quasi-morphisms and the Poisson bracket, \emph{Pure and Appl. Math.
Quarterly} \textbf{3} (2007), 1037--1055.

\bibitem{EPZ-physics} Entov, M., Polterovich, L., Zapolsky, F.,
An "anti-Gleason" phenomenon and simultaneous measurements in
classical mechanics, \emph{Found. of Physics} \textbf{37} (2007),
1306--1316.

\bibitem{FOOO-book} Fukaya, K., Oh, Y.-G., Ohta, H., Ono, K.,
\emph{Lagrangian intersection Floer theory: anomaly and obstruction.
Parts I, II.} AMS, Providence, RI, International Press, Somerville,
MA, 2009.


\bibitem{FOOO-spectral-with-bulk} Fukaya, K., Oh, Y.-G., Ohta, H., Ono, K.,
Spectral invariants with bulk, quasimorphisms and Lagrangian Floer theory,
arXiv:1105.5123, 2011.


\bibitem{Galkin} Galkin, S., The conifold point, preprint, arXiv:1404.7388, 2014.

\bibitem{Gambaudo-Ghys} Gambaudo, J.-M., Ghys, E., Commutators and
diffeomorphisms of surfaces, \emph{Erg. Th. Dyn.
Sys.} \textbf{24} (2004), 1591--1617.

\bibitem{Givental1} Givental, A., The nonlinear Maslov index, in \emph{Geometry of low-dimensional manifolds, 2 (Durham, 1989), 35--43}. London Math. Soc. Lecture Note Ser.~151, Cambridge Univ. Press, Cambridge, 1990.

\bibitem{Givental2}  Givental, A., Nonlinear generalization of the Maslov index, in \emph{Theory of singularities and its applications, 71--103}. Adv. Soviet Math.~1, AMS, Providence, RI, 1990.

\bibitem{Gleason} Gleason, A.M., Measures on the
closed subspaces of a Hilbert space, \emph{J. Math. Mech.}
\textbf{6} (1957), 885--893.

\bibitem{Hastings} Hastings, M.B.,  Making almost commuting matrices commute,
\emph{Comm. Math. Phys.}  \textbf{291}  (2009),  321--345.


\bibitem{Kawasaki} Kawasaki, M., Superheavy subsets and noncontractible
Hamiltonian circle actions, preprint, 2013.

\bibitem{Khanevsky} Khanevsky, M.,
Hofer's metric on the space of diameters, \emph{J. Topol. and Analysis} \textbf{1} (2009), 407--416.

\bibitem{Kislev} Kislev, A., Compactly supported Hamiltonian loops with non-zero Calabi invariant, preprint, arXiv:1310.1555, 2013.

\bibitem{Knudsen-Advances} Knudsen, F.F., Topology and the construction of extreme quasi-measures,
\emph{Adv. Math.} \textbf{120} (1996), 302--321.



\bibitem{Lanzat} Lanzat, S., Quasi-morphisms and symplectic quasi-states for convex symplectic manifolds, \emph{Int. Math. Res. Not.} (2013), 5321-5365.

\bibitem{LeRoux} Le Roux, F., Six questions, a proposition and two pictures on Hofer distance for Hamiltonian diffeomorphisms on surfaces, in \emph{Symplectic topology and measure preserving dynamical systems, 33--40}. Contemp. Math.~512, AMS, Providence, RI, 2010.

\bibitem{Maydanskiy-Mirabelli} Maydanskiy, M., Mirabelli, B.P., Semisimplicity of the quantum cohomology for smooth Fano toric varieties associated with facet symmetric polytopes, \emph{Electron. Res. Announc. Math. Sci.} \textbf{18} (2011), 131--143.

\bibitem{McDuff-monodromy} McDuff, D., Monodromy in Hamiltonian Floer theory, \emph{Comment. Math. Helv.} \textbf{85} (2010), 95--133.

\bibitem{McDuff-probes} McDuff, D., Displacing Lagrangian toric fibers via probes, in \emph{Low-dimensional and symplectic topology, 131--160}. Proc. Sympos. Pure Math.~82, AMS, Providence, RI, 2011.

\bibitem{McD-Sal-intro} McDuff, D., Salamon, D., \emph{Introduction to symplectic topology}. Second edition. Oxford Univ. Press, New York, 1998.

\bibitem{McD-Sal-psh-book} McDuff, D., Salamon, D., \emph{$J$-holomorphic curves and symplectic topology}. Second edition. AMS Colloquium Publ.~52, AMS, Providence, RI, 2012.

\bibitem{Monzner-Vichery-Zapolsky} Monzner, A., Vichery, N., Zapolsky, F., Partial quasimorphisms and quasistates on cotangent bundles, and symplectic homogenization, \emph{J. Mod. Dyn.} \textbf{6} (2012), 205--249.

\bibitem{Oakley-Usher} Oakley, J., Usher, M., On certain Lagrangian submanifolds of $S^2\times S^2$ and $\C P^n$, preprint,  arXiv:1311.5152, 2013.

\bibitem{Oh1}  Oh, Y.-G.,  Symplectic topology as the
geometry of action functional I, \emph{J. Diff. Geom.}  \textbf{46}
(1997), 499--577.

\bibitem{Oh2} Oh, Y.-G.,  Symplectic topology as the
geometry of action functional II, \emph{Comm. Analysis Geom.} \textbf{7}
(1999), 1--55.

\bibitem{Oh-spectral} Oh, Y.-G.,
Construction of spectral invariants of Hamiltonian
paths on general symplectic manifolds, in \emph{The
breadth of symplectic and Poisson geometry, 525--570}.
Birkh\"auser, Boston, 2005.

\bibitem{Ostr-Calabi} Ostrover, Y., Calabi quasi-morphisms for some non-monotone symplectic manifolds,
\emph{Algebr. Geom. Topol.} \textbf{6} (2006), 405--434.

\bibitem{Ostrover-Tyomkin} Ostrover, Y., Tyomkin, I.,
On the quantum homology algebra of toric Fano manifolds, \emph{Selecta Math. (N.S.)} \textbf{15} (2009), 121--149.



\bibitem{P-quantum-unsharpness} Polterovich, L., Quantum unsharpness and symplectic rigidity, \emph{Lett. Math. Phys.} \textbf{102} (2012), 245--264.

\bibitem{P-quantum-noise} Polterovich, L., Symplectic geometry of quantum noise, preprint,  arXiv:1206.3707, 2012. To appear in \emph{Comm. Math. Phys.}

\bibitem{PR-book} Polterovich, L., Rosen., D., \emph{Function theory on symplectic manifolds}, book draft, 2014.

\bibitem{Py-AnnENS} Py, P., Quasi-morphismes et invariant de Calabi, \emph{Ann. Sci. \'Ecole Norm. Sup. (4)} \textbf{39} (2006), 177--195.

\bibitem{Py-CR-torus}  Py, P., Quasi-morphismes de Calabi et graphe de Reeb sur le tore, \emph{C. R. Math. Acad. Sci. Paris} \textbf{343} (2006), 323--328.

\bibitem{Schwarz} Schwarz, M.,
 On the action spectrum for closed symplectically aspherical
manifolds, \emph{Pacific J. Math.} \textbf{193} (2000), 419--461.

\bibitem{Seyfaddini} Seyfaddini, S., Descent and $C^0$-rigidity of spectral invariants on monotone symplectic manifolds, \emph{J. Topol. Analysis} \textbf{4} (2012), 481--498.

\bibitem{Seyfaddini-LagrHofer} Seyfaddini, S., Unboundedness of the Lagrangian Hofer distance in the Euclidean ball, preprint,  arXiv:1310.1057, 2013.

\bibitem{Shelukhin} Shelukhin, E., Remarks on invariants of Hamiltonian loops, \emph{J. Topol. Analysis} \textbf{2} (2010), 277--325.

\bibitem{Shelukhin-qmms} Shelukhin, E., The action homomorphism, quasimorphisms and moment maps on the space of compatible almost complex structures, preprint, arXiv:1105.5814, 2011.

\bibitem{Shtern}  Shtern, A.I., The Kazhdan-Milman problem for semisimple compact Lie groups, \emph{Russian Math. Surveys}
\textbf{62} (2007), 113--174.

\bibitem{Tamarkin} Tamarkin, D., Microlocal condition for non-displaceablility, preprint, arXiv:0809.1584, 2008.

\bibitem{Usher-spectral} Usher, M., Spectral numbers in Floer theories, \emph{Compositio Math.} \textbf{144} (2008), 1581--1592.

\bibitem{Usher-duality} Usher, M., Duality in filtered Floer-Novikov complexes, \emph{J. Topol. Analysis} \textbf{2} (2010), 233--258.

\bibitem{Usher-deformed} Usher, M., Deformed Hamiltonian Floer theory, capacity estimates and Calabi quasimorphisms, \emph{Geom. Topol.} \textbf{15} (2011), 1313--1417.

\bibitem{Viterbo-spectral} Viterbo, C.,
 Symplectic topology as the geometry of generating functions,
\emph{Math. Ann.} \textbf{292} (1992), 685--710.

\bibitem{von Neumann} Von Neumann, J., \emph{Mathematical foundations of quantum
mechanics}. Princeton University Press, Princeton, 1955.
(Translation of \emph{Mathematische Grundlagen der
Quantenmechanik}. Springer, Berlin, 1932.)

\bibitem{Wu} Wu, W., On an exotic Lagrangian torus in $\C P^2$, preprint, arXiv:1201.2446, 2012.

\end{thebibliography}

\end{document}